\documentclass{amsart}
\usepackage{verbatim}
\usepackage{amssymb}
\usepackage{amsbsy}
\usepackage{amscd}
\usepackage{amsmath}
\usepackage{amsthm}
\usepackage[mathscr]{eucal}

\newenvironment{NB}{
{\bf NB}. \footnotesize
}{}
\renewenvironment{NB}{
\comment
  }{\endcomment}

\theoremstyle{plain}
 \newtheorem{thm}{Theorem}[section]
 \newtheorem{lem}[thm]{Lemma}
 \newtheorem{prop}[thm]{Proposition}
 \newtheorem{cor}[thm]{Corollary}
 
\theoremstyle{definition}
 \newtheorem{defn}{Definition}[section]
\theoremstyle{remark}
 \newtheorem{rem}{Remark}[section]
 
 \newtheorem{claim}{Claim}[section]
\def\Bbb{\mathbb}
\def\frak{\mathfrak}
\def\cal{\mathcal}

\newcommand{ \Supp}{\operatorname{Supp}}
\newcommand{\Ext}{\operatorname{Ext}}
\newcommand{\Hom}{\operatorname{Hom}}

\newcommand{\im}{\operatorname{im}}

\newcommand{\rk}{\operatorname{rk}}

\newcommand{\NS}{\operatorname{NS}}
\newcommand{\coker}{\operatorname{coker}}
\newcommand{\Pic}{\operatorname{Pic}}
\newcommand{\ch}{\operatorname{ch}}
\newcommand{\td}{\operatorname{td}}

\newcommand{\Hilb}{\operatorname{Hilb}}
\newcommand{\Quot}{\operatorname{Quot}}

\newcommand{\Coh}{\operatorname{Coh}}
\newcommand{\Spec}{\operatorname{Spec}}
\newcommand{\Amp}{\operatorname{Amp}}

\newcommand{\id}{\operatorname{id}}

\newcommand{\Br}{\operatorname{Br}}
\newcommand{\tr}{\operatorname{tr}}
\newcommand{\dslash}{/\!\!/} % for algebro-geometric quotient (double
\font\b=cmr10 scaled \magstep5
\def\bigzerou{\smash{\lower1.7ex\hbox{\b 0}}}

\numberwithin{equation}{section}
\begin{document}

\title
{
Moduli spaces of twisted sheaves on a projective variety}
%\dedicatory{Dedicated to Masaki Maruyama on the occation 
%of his 60th birthday}
\author{K\={o}ta Yoshioka
% Department of Mathematics\\
%Faculty of Science, Kobe University
}
 \address{Department of Mathematics, Faculty of Science,
Kobe University,
Kobe, 657, Japan}
\email{yoshioka@math.kobe-u.ac.jp}
 
\maketitle

\section{Introduction}
%In this note, we construct the moduli space of twisted sheaves
%by modifying \cite{Y:11}.
%We also generalize some results in \cite{Y:7} to this case.

Let $X$ be a smooth projective variety over ${\Bbb C}$. 
Let $\alpha:=\{\alpha_{ijk} 
\in H^0(U_i \cap U_j \cap U_k,{\cal O}_X^{\times})\}$ be a 2-cocycle
representing a torsion class $[\alpha] \in H^2(X,{\cal O}_X^{\times})$.
An $\alpha$-twisted sheaf 
$E:=\{(E_i,\varphi_{ij})\}$ is a collection of sheaves
$E_i$ on $U_i$ and isomorphisms $\varphi_{ij}:
E_{i|U_i \cap U_j} \to E_{j|U_i \cap U_j}$
such that $\varphi_{ii}=\id_{E_i}$, $\varphi_{ji}=\varphi_{ij}^{-1}$
and 
$\varphi_{ki} \circ \varphi_{jk} \circ \varphi_{ij} 
=\alpha_{ijk}\id_{E_i}$.
We assume that there is a locally free $\alpha$-twisted sheaf,
that is, $\alpha$ gives an element of the Brauer group $\Br(X)$.
A twisted sheaf naturally appears if we consider a non-fine moduli space
$M$ of the usual stable sheaves on $X$.
Indeed the transition functions of the local universal families 
satisfy the patching condition up to the multiplication by
constants and gives a twisted sheaf.
If the patching condition is satisfied, i.e.,
the moduli space $M$ is fine, than the universal family
defines an integral functor on the bounded 
derived categories of coherent sheaves
${\bf D}(M) \to {\bf D}(X)$.
Assume that $X$ is a $K3$ surface and $\dim M=\dim X$.
Than
Mukai, Orlov and Bridgeland showed that
the integral functor is the Fourier-Mukai functor, i.e.,
it is an equivalence of the categories.
In his thesis \cite{C:2}, 
C\u ald\u araru studied the category of twisted sheaves and 
its bounded derived category.
In particular, he generalized Mukai, Orlov and Bridgeland's
results on the Fourier-Mukai transforms to non-fine moduli spaces
on a $K3$ surface.
For the usual derived category, Orlov \cite{Or:1} showed that
the equivalence class is described in terms of the Hodge
structure of the Mukai lattice.
C\u ald\u araru conjectured that a similar result also holds
for the derived category of twisted sheaves.
Recently this conjecture was modified and proved by 
Huybrechts and Stellari, if $\rho(X) \geq 12$ in \cite{H-S:2}.
As is well-known,
twisted sheaves also appear if we consider a projective bundle over $X$.

%This means that $\alpha$ gives a torsion class of 
%$H^2(X,{\cal O}_X^{\times})$.
%Conversely for a 2-cocycle $\alpha$ which gives a torsion class 
%of $H^2(X,{\cal O}_X^{\times})$,
%Grothendieck asked that there is a locally free $\alpha$-sheaf. 
%If $X$ is a projective surface or a $K3$ surface,
%$\Br(X)$ is the torsion part of $H^2(X,{\cal O}_X^{\times})$.
%this question is solved affirmatively.

In this paper, we define a notion of the stability for 
a twisted sheaf and construct
the moduli space of stable twisted sheaves on $X$.
We also construct a projective compactification of
the moduli space by adding the $S$-equivalence classes
of semi-stable twisted sheaves. 
In particular if $H^1(X,{\cal O}_X)=0$ (e.g. $X$ is a $K3$ surface), 
then the moduli space of locally free twisted sheaves is  
the moduli space of projective bundles over $X$.
Thus we compactify the moduli space of projective bundles
by using twisted sheaves. 
The idea of the construction is as follows.
We consider a twisted sheaf as a usual sheaf on the Brauer-Severi
variety. 
Instead of using the Hilbert polynomial associated to
an ample line bundle on the Brauer-Severi variety,
we use the Hilbert polynomial associated to
a line bundle coming from $X$ in order to define the stability. 
Then the construction is a modification of Simpson's construction
of the moduli space of usual sheaves (cf. \cite{Y:11}). 
M. Lieblich informed us that our stability condition coincides with
Simpson's stability for modules over the associated Azumaya algebra
via Morita equivalence. Hence the construction also follows from
Simpson's moduli space 
\cite[Thm. 4.7]{S:1} and the valuative criterion for properness. 
  
In section \ref{sect:k3}, 
we consider the moduli space of twisted sheaves on a $K3$ surface.
We generalize known results on the moduli space of usual stable sheaves
to the moduli spaces of twisted stable sheaves (cf. \cite{Mu:3},
\cite{Y:7}).
In particular, we consider the non-emptyness, the deformation type 
and the weight 2 Hodge structure.
Then we can consider twisted version of the Fourier-Mukai transform
by using 2 dimensional moduli spaces, 
which is done in section \ref{sect:FM}.  
As an application of our results, 
Huybrechts and Stellari \cite{H-St:2}
prove C\u ald\u araru's conjecture
generally. 

Since our main example of twisted sheaves
are those on $K3$ surfaces or abelian surfaces, 
we consider twisted sheaves over ${\Bbb C}$.
But some of the results (e.g., subsection \ref{subsect:construction})
also hold over any field.
 
E. Markman and D. Huybrechts communicated to the author that
M. Lieblich independently constructed 
the moduli of twisted sheaves.
In his paper \cite{Li:1},
Lieblich developed a general theory of twisted sheaves
in terms of algebraic stack and constructed the moduli space
intrinsic way.
He also studied the moduli spaces of twisted sheaves on 
surfaces.
There are also some overlap with a paper by 
N. Hoffmann and U. Stuhler \cite{Ho-St:1}.
They also constructed the moduli space of rank 1 twisted sheaves
and studied the symplectic structure of the moduli space.  
\begin{NB}
\subsection{}

\begin{prop}\cite[Prop. 2.2]{D:1}
$H^1(PGL(r)) \to H^2(X,\mu_r)$ is surjective.
\end{prop}

By the exact and commutative diagram
\begin{equation}
\begin{CD}
@. 1 @.  1 @. @.\\
@. @VVV @VVV @. @.\\
1 @>>> \mu_r @>>> SL(r) @>>> PSL(r) @>>> 1\\
@. @VVV @VVV @| @. \\
1 @>>> {\cal O}_X^{\times} @>>> GL(r) @>>> PGL(r) @>>> 1\\
@. @V{(\;\;)^r}VV @VV{\det}V @. @.\\
@. {\cal O}_X^{\times} @= {\cal O}_X^{\times} @. @.\\
@. @VVV @VVV @. @.\\
@. 1 @.  1 @. @.,
\end{CD}
\end{equation}
we have an exact and commutative diagram

\begin{equation}
\begin{CD}
@. @. @. H^1(X,{\cal O}_X^{\times})\\
@. @. @. @VV{c_1 \mod r}V \\
H^1(X,\mu_r) @>>> H^1(X,SL(r)) @>>> H^1(PGL(r)) @>{\delta'}>> H^2(X,\mu_r)\\
@VVV @VVV @| @VV{o}V \\
H^1(X,{\cal O}_X^{\times}) 
@>>> H^1(X,GL(r)) @>>> H^1(PGL(r)) @>{\delta}>> H^2(X,{\cal O}_X^{\times})
\end{CD}
\end{equation}

\begin{equation}
\begin{CD}
@. 0 @.  0 @. @.\\
@. @VVV @VVV @. @.\\
@. {\Bbb Z} @=  {\Bbb Z} @. @.\\
@. @V{r}VV @VV{r}V @. @.\\
0 @>>> {\Bbb Z} @>>> {\cal O}_X @>>> {\cal O}_X^{\times} @>>> 1\\
@. @VVV @VVV @| @. \\
1 @>>> \mu_r @>>> {\cal O}_X^{\times} @>{(\;\;)^r}>> 
{\cal O}_X^{\times} @>>> 1\\
@. @VVV @VVV @. @.\\
@. 1 @.  1 @. @.\\
\end{CD}
\end{equation}

\begin{equation}
\begin{CD}
H^1(X,{\cal O}_X^{\times}) @>{c_1}>> H^2(X,{\Bbb Z}) @>>>
H^2(X,{\cal O}_X)\\
@| @VVV @VVV \\
H^1(X,{\cal O}_X^{\times}) @>>> H^2(X,\mu_r) @>{\lambda}>>
H^2(X,{\cal O}_X^{\times})
\end{CD}
\end{equation}
Hence
\begin{equation}
o(H^2(X,\mu_r))=H^2(X,{\Bbb Z})/(\NS(X)+rH^2(X,{\Bbb Z})).
\end{equation}

\end{NB}

%{\bf Notation}
%For a ${\Bbb P}^{r-1}$-bundle $p:Y \to X$, 
%$[Y] \in H^1(X,PGL(r))$ denotes the corresponding cohomology class.

\section{Twisted sheaves}
Notation:
For a locally free sheaf $E$ on a variety $X$,
${\Bbb P}(E) \to X$ denotes the projective bundle
in the sense of Grothendieck, that is, 
${\Bbb P}(E)=\mathrm{Proj}(\bigoplus_{n=0}^{\infty} S^n(E))$.

\vspace{1pc}

Let $X$ be a smooth projective variety over ${\Bbb C}$. 
Let $\alpha:=\{\alpha_{ijk} 
\in H^0(U_i \cap U_j \cap U_k,{\cal O}_X^{\times})\}$ be a 2-cocycle
representing a torsion class $[\alpha] \in H^2(X,{\cal O}_X^{\times})$.
An $\alpha$-twisted sheaf 
$E:=\{(E_i,\varphi_{ij})\}$ is a collection of sheaves
$E_i$ on $U_i$ and isomorphisms $\varphi_{ij}:
E_{i|U_i \cap U_j} \to E_{j|U_i \cap U_j}$
such that $\varphi_{ii}=\id_{E_i}$, $\varphi_{ji}=\varphi_{ij}^{-1}$
and 
$\varphi_{ki} \circ \varphi_{jk} \circ \varphi_{ij} 
=\alpha_{ijk}\id_{E_i}$.
If all $E_i$ are coherent, then we say that $E$ is coherent. 
Let $\Coh(X,\alpha)$ be the category of coherent 
$\alpha$-twisted sheaves 
on $X$.

If $E_i$ are locally free for all $i$, 
then we can glue ${\Bbb P}(E_i^{\vee})$ together
and get a projective bundle $p:Y \to X$ with
$\delta([Y])=[\alpha]$, where $[Y] \in H^1(X,PGL(r))$ is the corresponding
cohomology class of $Y$ and
$\delta:H^1(X,PGL(r)) \to 
H^2(X,{\cal O}_X^{\times})$ is the connecting homomorphism induced by
the exact sequence
\begin{equation}
1 \to {\cal O}_X^{\times} \to GL(r) \to PGL(r) \to 1.
\end{equation}
Thus $[\alpha]$ belongs to the Brauer group $\Br(X)$.
If $X$ is a smooth projective surface, then
$\Br(X)$ coincides with 
the torsion part of $H^2(X,{\cal O}_X^{\times})$.
Let ${\cal O}_{{\Bbb P}(E_i^{\vee})}(\lambda_i)$ be the tautological
line bundle on ${\Bbb P}(E_i^{\vee})$.
Then, $\varphi_{ij}$ induces an isomorphism
$\widetilde{\varphi}_{ij}:
{\cal O}_{{\Bbb P}(E_i^{\vee})}(\lambda_i)_{|p^{-1}(U_i \cap U_j)} 
\to {\cal O}_{{\Bbb P}(E_j^{\vee})}(\lambda_j)_{|p^{-1}(U_i \cap U_j)}$.
${\cal L}(p^*(\alpha^{-1})):=\{({\cal O}_{{\Bbb P}(E_i^{\vee})}(\lambda_i),
\widetilde{\varphi}_{ij})\}$ is an $p^*(\alpha^{-1})$-twisted line bundle
on $Y$.

\subsection{Sheaves on a projective bundle}
In this subsection, we shall interpret twisted sheaves as
usual sheaves on a Brauer-Severi variety.
Let $p:Y \to X$ be a projective bundle. 
\begin{NB}
Since $-K_Y$ is relatively ample,
$Y$ is projective over $X$.
\end{NB}
Let $X=\cup_i U_i$ be an analytic open covering of $X$
such that $p^{-1}(U_i) \cong U_i \times {\Bbb P}^{r-1}$.
We set $Y_i:=p^{-1}(U_i)$.
We fix a collection of tautological line bundles
${\cal O}_{Y_i}(\lambda_i)$ on $Y_i$
and isomorphisms
$\phi_{ji}:{\cal O}_{Y_i \cap Y_j}(\lambda_j) \to
{\cal O}_{Y_i \cap Y_j}(\lambda_i)$.
We set $G_i:=p_*({\cal O}_{Y_i}(\lambda_i))^{\vee}$.
Then $G_i$ are vector bundles on $U_i$ and
$p^*(G_i)(\lambda_i)$ defines a vector bundle 
$G$ of rank $r$ on $Y$.
We have the Euler sequence
\begin{equation}
0 \to {\cal O}_Y \to G \to T_{Y/X} \to 0.
\end{equation}
Thus $G$ is a non-trivial extension of $T_{Y/X}$ by
${\cal O}_Y$.

\begin{lem}\label{lem:G}
$\Ext^1(T_{Y/X},{\cal O}_Y)={\Bbb C}$.
Thus $G$ is characterized as a 
non-trivial extension of $T_{Y/X}$ by
${\cal O}_Y$.
In particular, $G$ does not depend on the choice of the 
local trivialization of $p$.
\end{lem}
\begin{proof}
Since ${\bf R} p_* (G^{\vee})=0$,
the Euler sequence inplies that 
$\Ext^1(T_{Y/X},{\cal O}_Y) \cong H^0(Y,{\cal O}_Y) \cong {\Bbb C}$.
\end{proof}

\begin{defn}
For a projective bundle $p:Y \to X$,
let $\epsilon(Y)(:=G)$ be a vector bundle on $Y$ 
which is a non-trivial extension
\begin{equation}
0 \to {\cal O}_Y \to  \epsilon(Y)\to T_{Y/X} \to 0.
\end{equation}
\end{defn}

By the exact sequence
$0 \to \mu_r \to SL(r) \to PGL(r) \to 1$,
we have a connecting homomorphism
$\delta':H^1(X,PGL(r)) \to H^2(X,\mu_r)$. 
Let $o:H^2(X,\mu_r) \to H^2(X,{\cal O}_X^{\times})$
be the homomorphism induced by the inclusion
$\mu_r \hookrightarrow {\cal O}_X^{\times}$.
Then we have $\delta=o \circ \delta'$.
\begin{defn}\label{defn:w(Y)}
For a ${\Bbb P}^{r-1}$-bundle $p:Y \to X$ corresponding
to $[Y] \in H^1(X,PGL(r))$,
we set $w(Y):=\delta'([Y]) \in H^2(X,\mu_r)$.
\end{defn}

\begin{lem}[\cite{C:1},\cite{H-S:1}]
If $p:Y \to X$ is a ${\Bbb P}^{r-1}$-bundle associated to a vector bundle
$E$ on $X$, i.e., $Y={\Bbb P}(E^{\vee})$, then
$w(Y)=[c_1(E) \mod r]$.
\end{lem}
 
\begin{NB}
\begin{lem}
\begin{enumerate}
\item
$\Ext_p^i(T_{Y/X},{\cal O}_Y)=0$, $i \ne 1$ and
$\Ext_p^1(T_{Y/X},{\cal O}_Y) \otimes k(x) \cong 
\Ext_p^1(T_{Y/X|p^{-1}(x)},{\cal O}_{p^{-1}(x)})$. 
\item
$\Ext^1(T_{Y/X},{\cal O}_Y) \cong H^0(X,\Ext_p^1(T_{Y/X},{\cal O}_Y)) 
\cong {\Bbb C}$
\end{enumerate}
\end{lem}
We have
$K_{Y/X} \cong \det G^{\vee}$.

\end{NB}

\begin{NB}
If $w(Y)=[D \mod r] \in H^2(X,\mu_r)$ for $D \in \Pic(X)$, then
$[c_1(G) \mod r]=w(Y \times_X Y)=p^*(w(Y)) =p^*(D) \in
H^2(Y,\mu_r)$. Hence there is a line bundle $L$ on $X$
such that $c_1(G)-p^*(D)=rc_1(L)$. 
Then $L_{|Y_i} \cong {\cal O}_{Y_i}(1)$.
Hence $Y \cong {\Bbb P}(E^{\vee})$ for $E:=p_*(L)$ and
$G \cong E^{\vee} \otimes L$.
\end{NB}

\begin{lem}\label{lem:w(Y)}
$[c_1(G) \mod r]=p^*(w(Y)) \in H^2(Y,\mu_r)$.
\end{lem}

\begin{proof}
There is a line bundle $L$ on $Y \times_X Y$ such that
$L_{|Y_i \times_{U_i} Y_i} \cong
p_1^*({\cal O}_{Y_i}(-\lambda_i)) \otimes
p_2^*({\cal O}_{Y_i}(\lambda_i))$, where
$p_i:Y \times_X Y \to Y$, $i=1,2$ are $i$-th projections.
By the definition of $G$, 
$p_{1*} (L) \cong G^{\vee}$.
Hence $p_1: Y \times_X Y \to Y$ is the projective bundle
${\Bbb P}(G^{\vee}) \to Y$. Then we get
$-[c_1(G^{\vee}) \mod r]=w(Y \times_X Y)=p^*(w(Y))$.
\end{proof}

 \begin{lem}
 Let $p:Y \to X$ be a ${\Bbb P}^{r-1}$-bundle.
 Then the following conditions are equivalent.
 \begin{enumerate}
 \item[(1)]
 $Y={\Bbb P}(E^{\vee})$ for a vector bundle on $X$.
 \item[(2)]
 $w(Y) \in \NS(X) \otimes \mu_r$.
 \item[(3)]
 There is a line bundle $L$ on $Y$ such that
 $L_{|p^{-1}(x)} \cong {\cal O}_{p^{-1}(x)}(1)$.
 \end{enumerate}
 \end{lem}
 
 \begin{proof}
 $(2) \Rightarrow (3)$:
 If $w(Y)=[D \mod r]$, $D \in \NS(X)$, then 
 $c_1(\epsilon(Y))-p^*(D) \equiv 0 \mod r$. We take a line bundle
 $L$ on $Y$ with 
 $c_1(\epsilon(Y))-p^*(D)=rc_1(L)$.
 $(3) \Rightarrow (1)$: We set $E^{\vee}:=p_*(L)$.
 Then $Y={\Bbb P}(E^{\vee})$.
 \end{proof}

\begin{defn}
$\Coh(X,Y)$ is a subcategory of $\Coh(Y)$ such that 
$E \in \Coh(X,Y)$ if and only if
\begin{equation}
E_{|Y_i} \cong p^*(E_i) \otimes 
{\cal O}_{Y_i}(\lambda_i)
\end{equation}
 for $E_i \in \Coh(U_i)$.
For simplicity, we call $E \in \Coh(X,Y)$ a $Y$-sheaf.
\end{defn}
By this definition, $\{(U_i,E_i)\}$ gives a twisted sheaf on $X$.
Thus we have an equivalence
\begin{equation}\label{eq:Y=X}
\begin{matrix}
\Lambda^{{\cal L}(p^*(\alpha^{-1}))}:\Coh(X,Y)& \cong & \Coh(X,\alpha)\\
E & \mapsto & p_*(E \otimes L^{\vee}),
\end{matrix}
\end{equation}
where ${\cal L}(p^*(\alpha^{-1})):=\{({\cal O}_{Y_i}(\lambda_i),\phi_{ij})\}$
is a twisted line bundle on $Y$ and
$\alpha_{ijk}^{-1}\id_{{\cal O}_{Y_i}(\lambda_i)}
=\phi_{ki} \circ \phi_{jk} \circ \phi_{ij}$.

We have the following relations:
\begin{equation}
\begin{split}
p_*(G^{\vee} \otimes E)_{|U_i}=&
p_*(p^*(G_i^{\vee}) \otimes {\cal O}_{Y_i}(-\lambda_i) \otimes
p^*(E_i)  \otimes {\cal O}_{Y_i}(\lambda_i))\\
=& p_* p^*(G_i^{\vee} \otimes E_i)=G_i^{\vee} \otimes E_i,
\end{split}
\end{equation}
\begin{equation}
\begin{split}
p_*(E)_{|U_i}=&
p_*(p^*(E_i)  \otimes {\cal O}_{Y_i}(\lambda_i))\\
=& E_i \otimes p_*({\cal O}_{Y_i}(\lambda_i))=G_i^{\vee} \otimes E_i.
\end{split}
\end{equation}

\begin{lem}
A coherent sheaf $E$ on $Y$ belongs to $\Coh(X,Y)$
if and only if 
$\phi:
p^* p_*(G^{\vee} \otimes E) \to  G^{\vee} \otimes E$ is an isomorphism.
In particular $E \in \Coh(X,Y)$ is an open condition.
\end{lem}

\begin{proof}
$\phi_{|Y_i}$ is the homomorphism
$p^* G_i^{\vee} \otimes p^* p_*(E(-\lambda_i)) \to
p^* G_i^{\vee} \otimes E(-\lambda_i)$.
Hence $\phi_{|Y_i}$ is an isomorphism if and only if
$p^* p_*(E(-\lambda_i)) \to E(-\lambda_i)$ is an isomorphism,
which is equivalent to $E \in \Coh(X,Y)$.
\end{proof}

\begin{lem}
Assume that $H^3(X,{\Bbb Z})_{\text{tor}}=0$.
Then $H^*(Y,{\Bbb Z}) \cong H^*(X,{\Bbb Z})[x]/(f(x))$,
where $f(x) \in H^*(X,{\Bbb Z})[x]$ is a monic polynomial
of degree $r$.
In particular, $H^2(X,{\Bbb Z}) \otimes \mu_{r'} \to
H^2(Y,{\Bbb Z}) \otimes \mu_{r'}$ is
injective for all $r'$.
\end{lem}
\begin{proof}
$R^2 p_*{\Bbb Z}$ is a local system of rank 1.
Since $c_1(K_{Y/X})$ is a section of this local system,
$R^2 p_*{\Bbb Z} \cong {\Bbb Z}$.
Let $h$ be the generator. Then $R^{2i} p_*{\Bbb Z} \cong {\Bbb Z}h^i$.
Since $H^3(X,{\Bbb Z})_{\text{tor}}=0$,
by the Leray spectral sequence, we get a surjective homomorphism
$H^2(Y,{\Bbb Z}) \to H^0(X,R^2 p_*{\Bbb Z})$.
Let $x \in H^2(Y,{\Bbb Z})$ be a lifting of $h$.
Then $x^i$ is a lifting of $h^i \in H^0(X,R^{2i} p_*{\Bbb Z})$.
Therefore the Leray-Hirsch theorem implies that
$H^*(Y,{\Bbb Z}) \cong H^*(X,{\Bbb Z})[x]/(f(x))$.
\end{proof}

\begin{NB}
If $H^3(X,{\Bbb Z})_{\text{tor}}=0$, then 
$H^3(Y,{\Bbb Z})_{\text{tor}}=0$.
$H^2(X,\mu_r) \cong H^2(X,{\Bbb Z}) \otimes \mu_r$.
\end{NB}

\begin{NB}
$Y \times_X Y \to Y$ is a projective bundle associated to
$G^{\vee}$:
$$
(1 \times p)^*(G^{\vee})_{|Y_i \times_{U_i} Y_i}= 
(1 \times p)^*({\cal O}_{Y_i}(-\lambda_i)) \otimes
 (p \times p)^*(G_i^{\vee})  \to
{\cal O}_{Y_i}(-\lambda_i) \boxtimes {\cal O}_{Y_i}(\lambda_i)
$$
is the tautological line bundle on $Y \times_X Y$.

$p':Y' \to X$ be a ${\Bbb P}^{r-1}$-bundle with a vector bundle $G'$ on $Y$.
Assume that $w(Y')=w(Y)$. Then
$c_1(G')-c_1(G) \mod r =0$ in $H^2(Y' \times_X Y,{\Bbb Z})$.
Then there is a line bundle $L$ on $Y' \times_X Y$
such that $rc_1(L)=c_1(G')-c_1(G)$.
For all fiber $l$ of $p' \times p:Y' \times_X Y \to X$, 
$L_{|l} \cong {\cal O}_{{\Bbb P}^{r-1}}(1) \boxtimes 
 {\cal O}_{{\Bbb P}^{r-1}}(-1)$.
Then we have an equivalence:
\begin{equation}
\begin{matrix} 
\Xi_{Y \to Y'}^L:\Coh(X,Y) & \to &\Coh(X,Y')\\ 
E & \mapsto & (p' \times 1)_*((1 \times p)^*(E) \otimes L).
\end{matrix}
\end{equation}

\end{NB}

\begin{lem}\label{lem:o(w(Y))}
Assume that $o(w(Y))=o(w(Y'))$.
\begin{enumerate}
\item
Then there is a line bundle $L$ on $Y' \times_X Y$
such that $L_{|{p'}^{-1}(x) \times p^{-1}(x)} \cong
{\cal O}_{{p'}^{-1}(x)}(1) \boxtimes {\cal O}_{p^{-1}(x)}(-1)$
for all $x \in X$.
If $L' \in \Pic(Y' \times_X Y)$ also satisfies the property, then
$L'=L \otimes q^*(P)$, $P \in \Pic(X)$, where $q:Y' \times_X Y \to X$
is the projection.  
\item
We have an equivalence
\begin{equation}
\begin{matrix} 
\Xi_{Y \to Y'}^L:&\Coh(X,Y) & \to &\Coh(X,Y')\\ 
& E & \mapsto & p_{Y'*}({p'}_{Y}^*(E) \otimes L),
\end{matrix}
\end{equation} 
where $p_{Y'}:Y' \times_X Y \to Y'$ and $p'_Y:Y' \times_X Y \to Y$
are projections.
\end{enumerate}
\end{lem}

\begin{NB}
We note that
$[r'c_1(G) \mod rr'], [rc_1(G') \mod rr']
 \in H^2(Y' \times_X Y,\mu_{rr'})$
belongs to $H^2(X,\mu_{rr'})$.
Since we have a commutative diagram
\begin{equation}
\begin{CD}
H^2(X,\mu_r) @>{o}>> H^2(X,{\cal O}_X^{\times})\\
@V{r'}VV @| \\
H^2(X,\mu_{rr'}) @>{o}>> H^2(X,{\cal O}_X^{\times}),
\end{CD}
\end{equation}
$o([r'c_1(G) \mod rr'])=o(w(Y))$.
Hence $[rc_1(G')-r'c_1(G) \mod rr'] \in H^2(X,\mu_{rr'})$
belongs to $\Pic(X) \otimes \mu_{rr'}$.
Therefore there is a line bundle $L$ on $Y' \times_X Y$
such that $rc_1(G')-r'c_1(G)=rr'c_1(L)+c_1(N)$, 
$N \in \Pic(X)$.
\end{NB}

\begin{rem}
We also see that
$E$ is a $Y$-sheaf if and only if
${p'}_{Y}^*(E) \otimes L \cong p_{Y'}^*(E')$ for a sheaf $E'$ on $Y'$.
\end{rem}

\begin{defn}\label{defn:w(E)}
Assume that $H^3(X,{\Bbb Z})_{tor}=0$.
For a $Y$-sheaf $E$ of rank $r'$,
$[c_1(E) \mod r'] \in H^2(Y,\mu_{r'})$ belongs to
$p^*( H^2(X,\mu_{r'}))$.
We set $w(E):=(p^*)^{-1}([c_1(E) \mod r']) \in H^2(X,\mu_{r'})$.
\end{defn}
By Lemmas \ref{lem:w(Y)} and \ref{lem:o(w(Y))}, we see that
\begin{lem}\label{lem:w(E)}
\begin{enumerate}
\item
By the functor $\Xi_{Y \to Y'}^L$ in Lemma \ref{lem:o(w(Y))},
$w(\Xi_{Y \to Y'}^L(E))=w(E)$ for $E \in \Coh(X,Y)$.
\item
$w(\epsilon(Y))=w(Y)$.
\end{enumerate}
\end{lem}

\begin{NB}
Let $\alpha:=\{\alpha_{ijk} \in 
H^0(U_i \cap U_j \cap U_k,{\cal O}_X^{\times})
\}$ and 
$\alpha':=\{\alpha'_{ijk}\in 
H^0(U_i \cap U_j \cap U_k,{\cal O}_X^{\times})\}$
are 2-cocyles such that they defines the same element of
the Brauer group:
$\alpha_{ijk}f_{ij}f_{jk}f_{ki}=\alpha_{ijk}'$.
Then we have an equivalence
\begin{equation}
\begin{matrix}
\Coh(X,\alpha)& \cong & \Coh(X,\alpha')\\
E & \mapsto & E \otimes {\cal O}_X^{\alpha^{-1} \alpha'}
\end{matrix}
\end{equation}
where ${\cal O}_X^{\alpha^{-1} \alpha'}$ is a $\alpha^{-1} \alpha'$-twisted
line bundle with ${\cal O}_X^{\alpha^{-1} \alpha'}=
\{({\cal O}_{U_i},f_{ij})\}$.

If there is a $\alpha$-twisted line bundle $L$, then
\begin{equation}
\begin{matrix}
\Coh(X)& \cong &\Coh(X,\alpha)\\
E & \mapsto & E \otimes L
\end{matrix}
\end{equation}

For a class $[\alpha] \in H^2(X,{\cal O}_X^{\times})$, 
we fix a representative $\alpha$ of $[\alpha]$.
Let $E=\{(E_i,\phi_{ij})\}$
and $E'=\{(E_i',\phi_{ij}')\}$ be $\alpha$-twisted vector bundles
such that ${\Bbb P}(E) \cong {\Bbb P}(E')$.
Then there is a collection of $\psi_i:E_i \to E_i'$
which fits in the commutative diagram:
\begin{equation}
\begin{CD}
E_{i|U_i \cap U_j} @>{\psi_i}>> E_{i|U_i \cap U_j}'\\
@V{\lambda_{ij}\phi_{ij}}VV @VV{\phi_{ij}'}V\\
E_{j|U_i \cap U_j} @>{\psi_j}>> E_{j|U_i \cap U_j}'
\end{CD}
\end{equation}
where $\lambda_{ij} \in H^0(U_i \cap U_j,{\cal O}_X^{\times})$.
Then we see that $\lambda_{ij}$ is a 1-cocycle. 
Let $L$ be a line bundle on $X$ accosiated to $\{\lambda_{ij}\}$.
Then $E' \cong E \otimes L$.
Therefore the moduli of projective bundles is the quotient of
the moduli of twisted vector bundles by the action of $\Pic(X)$.
\end{NB}

\section{Moduli of twisted sheaves}\label{sect:moduli}

\subsection{Definition of the stability}\label{subsect:defn}

Let $(X,{\cal O}_X(1))$ be a pair of
a projective scheme $X$ and an ample line bundle ${\cal O}_X(1)$
on $X$.
Let $p:Y \to X$ be a projective bundle over $X$.
\begin{defn}
A $Y$-sheaf $E$ is of dimension $d$, 
if $p_*(E)$ is of dimension $d$.
\end{defn}

For a coherent sheaf $F$ of dimension $d$ on $X$, 
we define $a_i(F) \in {\Bbb Z}$ by
the coefficient of the Hilbert polynomial of 
$F$: 
\begin{equation}
\chi(F(m))=
\sum_{i=0}^d a_i(F) \binom{m+i}{i}.
\end{equation}
Let $G$ be a locally free $Y$-sheaf.
For a $Y$-sheaf $E$ of dimension $d$, 
we set $a_i^{G}(E):=a_i(p_*(G^{\vee} \otimes E))$.
Thus we have
\begin{equation}
\chi(G,E \otimes p^*{\cal O}_X(m))=
\chi(p_*(G^{\vee} \otimes E)(m))=\sum_{i=0}^d a_i^{G}(E) \binom{m+i}{i}.
\end{equation}

\begin{defn}
Let $E$ be $Y$-sheaf of dimension $d$. Then
$E$ is ($G$-twisted) semi-stable (with respect to ${\cal O}_X(1)$),
if
\begin{enumerate}
\item
$E$ is of pure dimension $d$,
\item
\begin{equation}\label{eq:def}
\frac{\chi(p_*(G^{\vee} \otimes F)(m))}{a_d^{G}(F)} \leq 
\frac{\chi(p_*(G^{\vee} \otimes E)(m))}{a_d^{G}(E)},
m \gg 0
\end{equation}
for all subsheaf $F \ne 0$ of $E$.
\end{enumerate}
If the inequality in \eqref{eq:def} is strict for all
proper subsheaf $F \ne 0$ of $E$, then $E$ is ($G$-twisted) stable
with respect to ${\cal O}_X(1)$.
\end{defn}

\begin{thm}\label{thm:twisted}
Let $p:Y \to X$ be a projective bundle.
There is a coarse moduli scheme $\overline{M}_{X/{\Bbb C}}^{h}$
parametrizing $S$-equivalence classes of
$G$-twisted semi-stable $Y$-sheaves $E$ with the
$G$-twisted Hilbert polynomial $h$.
$\overline{M}_{X/{\Bbb C}}^{h}$ is a projective scheme.
\end{thm} 

\begin{rem}
The construction also works for a projective bundle
$Y \to X$ over any field and also for a family of projective bundles, 
by the fundamental work of Langer \cite{Langer:1}. 
\end{rem}
\begin{lem}\label{lem:choice}
Let $p':Y' \to X$ be a projective bundle with $o(w(Y'))=o(w(Y))$ and
$\Xi_{Y \to Y'}^L$ the correspondence in Lemma \ref{lem:o(w(Y))}.
Then a $Y$-sheaf $E$ is $G$-twisted semi-stable if and only if
$\Xi_{Y \to Y'}^L(E) \in \Coh(X,Y')$ is 
$\Xi_{Y  \to Y'}^L(G)$-twisted semi-stable.
In particular, we have an isomorphism of the corresponding moduli
spaces.
\end{lem}
Indeed, 
since $\Xi_{Y \times S \to Y'\times S}^{L \boxtimes {\cal O}_S}(*)_s=
\Xi_{Y \to Y'}^L(* \otimes k(s))$,
if we have a flat family of $Y$-sheaves 
$\{{\cal E}_s \}_{s \in S}$,
${\cal E} \in \Coh(Y \times S)$, then 
$\{{\cal E}_s'\}_{s \in S}$ is also a flat 
family of $Y'$-sheaves, where
${\cal E}':=
\Xi_{Y \times S \to Y'\times S}^{L \boxtimes {\cal O}_S}({\cal E})$.

\begin{rem}\label{rem:G}
For a locally free $Y$-sheaf $G$, we have a projective bundle
$Y' \to X$ with $\epsilon(Y')=\Xi_{Y \to Y'}^L(G)$.
Hence it is sufficient to study the $\epsilon(Y)$-twisted semi-stability. 
\end{rem}
\begin{rem}
This definition is the same as in \cite{C:1}.
If $Y={\Bbb P}(G^{\vee})$ for a vector bundle $G$ on $X$,
then $\Coh(X,Y)$ is equivalent to $\Coh(X)$ and
$G$-twisted stability is nothing but the 
twisted semi-stability in \cite{Y:11}.
\end{rem}

\begin{defn}
Let $\lambda$ be a rational number.
Let $E$ be a $Y$-sheaf of dimension $d$. Then
$E$ is of type $\lambda$
with respect to the $G$-twisted semi-stability,
if
\begin{enumerate}
\item
$E$ is of pure dimension $d$,
\item
\begin{equation}
\frac{a_{d-1}^{G}(F)}{a_d^{G}(F)} \leq 
\frac{a_{d-1}^{G}(E)}{a_d^{G}(E)}
+\lambda
\end{equation}
for all subsheaf $F$ of $E$.
\end{enumerate}
If $\lambda=0$, then $E$ is 
$\mu$-semi-stable.
\end{defn}

%\begin{defn}
%Let $H$ be an ample line bundle on $X$.
%$E \in \Coh(X,Y)$ is $Y$-twisted semi-stable with respect to $H$,
%if $E$ is torsion free and
%\begin{equation}
%\frac{\chi(p_*(F)(nH))}{\rk p_*(F)} \leq 
%\frac{\chi(p_*(G^{\vee} \otimes E)(nH))}{\rk p_*(G^{\vee} \otimes E)}\;(n \gg 0)
%\end{equation}
%for all subsheaf $F \in \Coh(X,Y)$ of $E$.
%\end{defn}

\subsection{Construction of the moduli space}\label{subsect:construction}
From now on, we assume that $G=\epsilon(Y)$ (cf. Remark \ref{rem:G}).
Let $P(x)$ be a numerical polynomial.
We shall construct the moduli space of $G$-twisted semi-stable $Y$-sheaves
$E$ with $\chi(p_*(G^{\vee} \otimes E)(n))=P(n)$.

\subsubsection{Boundedness}
Let $E$ be a $Y$-sheaf.
Then 
\begin{equation}
p^* p_*(G^{\vee} \otimes E) \otimes G \to E
\end{equation}
is surjective.
Indeed
$p^* p_*(G^{\vee} \otimes E) \to G^{\vee} \otimes E$
is an isomorphism and
$G \otimes G^{\vee} \to {\cal O}_{Y}$
is surjective.

We take a surjective homomorphism
${\cal O}_X(-n_G )^{\oplus N} \to p_*(G^{\vee} \otimes G)$, $n_G \gg 0$.
Then we have a surjective homomorphism
$p^*({\cal O}_X(-n_G ))^{\oplus N} \to G^{\vee} \otimes G$.

\begin{lem}\label{lem:bdd1}
Let $E$ be a $Y$-sheaf of pure dimension $d$.
If 
\begin{equation}\label{eq:type-nu}
a_{d-1}^G(F) \geq a_d^G(F)
\left(\frac{a_{d-1}^G(E)}{a_d^G(E)}-\nu \right)
\end{equation}
for all quotient $E \to F$, then
$a_{d-1}(F') \geq a_d(F') 
\left(\frac{a_{d-1}^G(E)}{a_d^G(E)}-\nu-n_G \right)$
for all quotient $p_*(G^{\vee} \otimes E) \to F'$.
In particular 
\begin{equation}
S_\nu:=\{E \in \Coh(X,Y)| \text{ $E$ satisfies \eqref{eq:type-nu}
and $\chi(p_*(G^{\vee} \otimes E)(nH))=P(n)$} \}
\end{equation}
is bounded.
\end{lem}

\begin{proof}
Since $p^* p_*(G^{\vee} \otimes E) \cong G^{\vee} \otimes E$,
we have a surjective homomorphism
\begin{equation}
p^*({\cal O}_X(-n_G H))^{\oplus N} \otimes E \to
G \otimes p^* p_*(G^{\vee} \otimes E) \to G \otimes p^*(F').
\end{equation}
By our assumption, we get 
\begin{equation}
a_{d-1}(p_*(G^{\vee} \otimes G) \otimes F') \geq 
a_d( p_*(G^{\vee} \otimes G) \otimes F')
\left(\frac{a_{d-1}(p_*(G^{\vee} \otimes E))}
{a_d(p_*(G^{\vee} \otimes E))}-n_G-\nu \right).
\end{equation}
Since $a_{d-1}(p_*(G^{\vee} \otimes G) \otimes F')=
\rk(G)^2 a_{d-1}(F')$ and
$a_{d}(p_*(G^{\vee} \otimes G) \otimes F')=
\rk(G)^2 a_{d}(F')$,
we get our claim.
The boundedness of $S_\nu$ follows from the boundedness of
$\{p_*(G^{\vee} \otimes E)| E \in S_{\nu} \}$ and Lemma \ref{lem:bdd2} below.
\end{proof}   

\begin{lem}\label{lem:bdd2}
Let $S$ be a bounded subset of $\Coh(X)$.
Then $T:=\{E \in \Coh(X,Y)| p_*(G^{\vee} \otimes E) \in S \}$
is also bounded.
\end{lem}

\begin{proof}
For $E \in T$, 
we set $I(E):=\ker(p^* p_*(G^{\vee} \otimes E) \otimes G \to E)$.
We shall show that $T':=\{ I(E)| E \in T \}$ is bounded.
We note that $I(E) \in \Coh(X,Y)$ and we have an exact sequence
\begin{equation}
0 \to p_*(G^{\vee} \otimes I(E)) \to  p_*(G^{\vee} \otimes E) 
\otimes p_*(G \otimes G^{\vee})
 \to
p_*(G^{\vee} \otimes E) \to 0.
\end{equation}
Since $p_*(G^{\vee} \otimes E) \in S$,
$\{ p_*(G^{\vee} \otimes I(E))| E \in T \}$ is also bounded.
Since $p^* p_*(G^{\vee} \otimes I(E)) \otimes G \to I(E)$ is surjective
and $I(E)$ is a subsheaf of 
$p^* p_*(G^{\vee} \otimes E) \otimes G$,
$T'$ is bounded.
\end{proof}

\begin{cor}
Under the same assumption \eqref{eq:type-nu},
there is a rational number $\nu'$ which depends on $\nu$
such that
\begin{equation}
a_{d-1}(F') \leq a_d(F')
\left(\frac{a_{d-1}^G(E)}{a_d^G(E)}+\nu' \right)
\end{equation}
for a subsheaf $F' \subset p_*(G^{\vee} \otimes E)$.
\end{cor}

Combining this with Langer's important result
\cite[Cor. 3.4]{Langer:1}, we have the following 
\begin{lem}\label{lem:section}
Under the same assumption \eqref{eq:type-nu},
\begin{equation}
 \frac{h^0(G,E)}{a_d^{G}(E)} \leq 
 \left[\frac{1}{d!}
 \left(\frac{a_{d-1}^{G}(E)}{a_d^{G}(E)}+\nu'
 +c \right)^d \right]_+,
\end{equation}
where $c$ depends only on $(X,{\cal O}_X(1))$, $G$, $d$ and $a_d^G(E)$.
\end{lem}

\subsubsection{A quot-scheme}

Since $p_*(G^{\vee} \otimes E)(n)$, $n \gg 0$ is generated by global sections, 
\begin{equation}
H^0(G^{\vee} \otimes E \otimes p^*{\cal O}_X(n))
 \otimes G \to E\otimes p^*{\cal O}_X(n)
\end{equation}
is surjective.
Since $R^i p_*(G^{\vee} \otimes E)=0$ for $i>0$, we also see that
$H^i(E \otimes p^*{\cal O}_X(n))=0$, $i>0$ and $n \gg 0$.

We fix a sufficiently large integer $n_0$.
We set $N:=\chi(p_*(G^{\vee} \otimes E)(n_0 ))=P(n_0)$.
We set $V:={\Bbb C}^N$.
We consider the quot-scheme ${\frak Q}$
parametrizing all quotients
\begin{equation}
\phi:V \otimes G \to E
\end{equation}
such that $E \in \Coh(X,Y)$ and 
$\chi(p_*(G^{\vee} \otimes E)(n))=P(n_0+n)$.
By Lemma \ref{lem:bdd2}, ${\frak Q}$ is bounded, in particular,
it is a quasi-projective scheme. 
\begin{lem}
${\frak Q}$ is complete.
\end{lem}

\begin{proof}
We prove our claim by using the valuative criterion.
Let $R$ be a discrete valuation ring
and $K$ the quotient field of $R$.
Let $\phi:V_R \otimes G \to {\cal E}$
be a $R$-flat family of quotients
such that ${\cal E} \otimes_R K \in \Coh(X,Y)$,
where $V_R:=V \otimes_{\Bbb C} R$.
We set ${\cal I}:=\ker \phi$.
We have an exact and commutative diagram:
\begin{equation}
\begin{CD}
0 @>>> p^* p_*({\cal I} \otimes G^{\vee}) @>>>
V_R \otimes G \otimes G^{\vee} @>>> p^* p_*({\cal E}\otimes G^{\vee})
@>>> 0\\
@. @VVV @| @VV{\psi}V \\
0 @>>> {\cal I} \otimes G^{\vee} @>>>
V_R \otimes G \otimes G^{\vee} @>>> {\cal E}\otimes G^{\vee}
@>>> 0
\end{CD}
\end{equation}
We shall show that $\psi$ is an isomorphism. 
Obviously $\psi$ is surjective.
Since ${\cal E}$ is $R$-flat, ${\cal E}$ has no $R$-torsion,
which implies that 
$p^* p_*({\cal E}\otimes G^{\vee})$ is a torsion free $R$-module.
Hence $\ker \psi$ is also torsion free.
On the other hand, our choice of ${\cal E}$ implies that
$\psi \otimes K$ is an isomorphism.
Therefore $\ker \psi=0$. 
\end{proof}

Since $\ker \phi \in \Coh(X,Y)$,
we have a surjective homomorphism
\begin{equation}
V \otimes \Hom(G,G \otimes p^*{\cal O}_X(n))
\to \Hom(G,E \otimes p^*{\cal O}_X(n))
\end{equation}
for $n \gg0$.
Thus we can embed ${\frak Q}$ as a subscheme of 
an Grassmann variety $Gr(V \otimes W,P(n_0+n))$,
where $W=\Hom(G,G \otimes p^*{\cal O}_X(n))$.
Since all semi-stable $Y$-sheaf are pure,
we may replace ${\frak Q}$ by the closure of the open subset
parametrizing pure quotient $Y$-sheaves. 
The same arguments in \cite{Y:11} imply that  
${\frak Q}\dslash GL(V)$ is the moduli space of $G$-twisted semi-stable sheaves.
The details are left to the reader.
For the proof, we also use the following. 

Let $(R,{\frak m})$ be a discrete valuation ring $R$ and the
maximal ideal ${\frak m}$.
Let $K$ be the fractional field and $k$ the residue field.
Let ${\cal E}$ be a $R$-flat family of $Y \otimes R$-sheaves
such that ${\cal E} \otimes_R K$
is pure.
\begin{lem}\label{lem:valuative}
There is a $R$-flat family
of coherent $Y \otimes R$-sheaves ${\cal F}$  and a homomorphism
$\psi:{\cal E} \to {\cal F}$ such that
${\cal F} \otimes_R k$ is pure,
$\psi_K$ is an isomorphism and
$\psi_k$ is an isomorphic at generic points of $\Supp({\cal F} \otimes_R k)$.
\end{lem}
By using \cite[Lem. 1.17]{S:1} or \cite[Prop. 4.4.2]{H-L:1},
we first construct ${\cal F}$ as a usual family of sheaves.
Then the very construction of it, ${\cal F}$ becomes a 
$Y \otimes R$-sheaf.

\begin{NB}
It also follows by using elementary transformations 
along $Y \otimes_R k$-sheaves on $Y \otimes_R k$. 
\end{NB}

\subsection{A family of $Y$-sheaves on a projective bundle over  
${M}_{X/{\Bbb C}}^{h}$}\label{subsubsect:family}

Assume that ${\frak Q}^{ss}$ consists of stable points.
Then ${\frak Q}^{ss} \to \overline{M}_{X/{\Bbb C}}^{h}$ is a 
principal $PGL(N)$-bundle.
For a scheme $S$, $f_S:Y \times S \to S$ denotes the projection.
Let ${\cal Q}$ be the universal quotient sheaf on
$Y \times {\frak Q}^{ss}$.
$V:=\Hom_{f_{ {\frak Q}^{ss}}}
(G \boxtimes {\cal O}_{{\frak Q}^{ss}},{\cal Q})$ is a locally free sheaf
on ${\frak Q}^{ss}$.
We consider the projective bundle 
${\frak q}:{\Bbb P}(V) \to {\frak Q}^{ss}$.
Since ${\cal Q}$ is $GL(N)$-linearized, $V$ is also $GL(N)$-linearized.
Then we have a quotient 
$\psi:{\Bbb P}(V) \to {\Bbb P}(V)/PGL(N)$ with the commutative diagram:
\begin{equation}
\begin{CD}
{\Bbb P}(V) @>{\frak q}>> {\frak Q}^{ss}\\
@VVV @VVV\\
\widetilde{\overline{M}_{X/{\Bbb C}}^{h}}:={\Bbb P}(V)/PGL(N) 
@>{q}>> \overline{M}_{X/{\Bbb C}}^{h}
\end{CD}
\end{equation}
Since $(1_Y \times {\frak q})^*({\cal Q}) \otimes 
f_{{\Bbb P}(V)}^*({\cal O}_{{\Bbb P}(V)}(-1))$ 
is $PGL(N)$-linearlized, we have a family of
$G$-twisted stable $Y$-sheaves ${\cal E}$ on 
$Y \times \widetilde{\overline{M}_{X/{\Bbb C}}^{h}}$
with $(1_Y \times \psi)^*({\cal E})=(1_Y \times {\frak q})^*({\cal Q}) \otimes 
f_{{\Bbb P}(V)}^*({\cal O}_{{\Bbb P}(V)}(-1))$.
Hence ${\cal E}^{\vee} \in \Coh(Y \times \overline{M}_{X/{\Bbb C}}^{h},
Y \times \widetilde{\overline{M}_{X/{\Bbb C}}^{h}})$
(if ${\cal E}$ is locally free).
Let $W$ be a locally free sheaf on 
$\widetilde{\overline{M}_{X/{\Bbb C}}^{h}}$
such that $\psi^*(W)={\frak q}^*(V)(-1)$.
Then we also have 
$W^{\vee}=\epsilon(\widetilde{\overline{M}_{X/{\Bbb C}}^{h}})
\in \Coh(\overline{M}_{X/{\Bbb C}}^{h},
\widetilde{\overline{M}_{X/{\Bbb C}}^{h}})$ and
${\cal E} \otimes f_{\widetilde{\overline{M}_{X/{\Bbb C}}^{h}}}^*(W^{\vee})$
descends to a sheaf on $Y \times \overline{M}_{X/{\Bbb C}}^{h}$.
 
\begin{rem}
There is also a family of $G$-twisted stable $Y$-sheaves ${\cal E}'$
on $Y \times {\Bbb P}(V^{\vee})/PGL(N)$
such that 
${\cal E}' \in \Coh(Y \times \overline{M}_{X/{\Bbb C}}^{h},
Y \times  {\Bbb P}(V^{\vee})/PGL(N))$.
\end{rem}

\begin{NB}

\subsection{Construction of moduli spaces}

From now on, we denote $a_{i}^{G}(*)$ by $a_{i}(*)$. 
For a $G$-twisted coherent sheaf $E \in \Coh(X,Y)$,
we denote $E \otimes p^* {\cal O}_X(n)$ by $E(n)$.
Let $E$ be a purely $d$-dimensional $G$-twisted coherent sheaf such that 
$E$ is of type $\lambda$ with respect to 
the $G$-twisted semi-stability.
Let $E \to E''$ be a quotient sheaf such that
$E''$ is of pure dimension $d$ and
\begin{equation}\label{eq:a3}
{a_d(E'')}\frac{\chi(G,E(n))}{a_d(E)} \geq 
{\chi(G,E''(n))}.
\end{equation}
Since the set of $E$ is bounded,
by Grothendieck's boundedness theorem, the set of such quotients
$E''$ is bounded.
Hence there is an integer $m({\lambda})$ which depends on $h$
and $\lambda$ such that, for $m \geq m({\lambda})$
and the kernel $E'$ of $E \to E''$ which satisfies \eqref{eq:a3},
\begin{enumerate}
 \item[($\flat 1$)]
  $\Hom(G,E'(m)) \otimes G \to E'(m)$ is surjective and
 \item[($\flat 2$)]
  $\Ext^i(G,E'(m))=0$, $i>0$.
\end{enumerate}
In particular, 
\begin{enumerate}
\item
$\Hom(G,E(m)) \otimes G \to E(m)$ is surjective and
\item 
 $\Ext^i(G,E(m))=0$, $i>0$.
\end{enumerate}

Let $V_m$ be a vector space of dimension $h(m)$.
Let ${\frak Q}:=\Quot_{V_m \otimes G/Y}^{h[m]}$ be the quot-scheme 
parametrizing all quotients
$V_m \otimes G \to F$ such that $F \in \Coh(X,Y)$ and 
the $G$-twisted Hilbert polynomial
of $F$ is $h[m]$, 
where $h[m](x)=h(m+x)$.
Let $V_m \otimes G \otimes {\cal O}_{{\frak Q}} \to \widetilde{E}(m)$
be the universal quotient sheaf on ${\frak Q} \times X$.
Let ${\frak Q}^{ss}$ be the open subscheme of ${\frak Q}$ consisting of
quotients $f:V_m \otimes G \to E(m)$
such that 
\begin{enumerate}
\item
a canonical map
$V_m \to \Hom(G,E(m))$ sending
$v \in V_m$ to $f(v \otimes *) \in \Hom(G,E(m))$ is an isomorphism and 
\item
$E$ is 
an $G$-twisted semi-stable sheaf.
\end{enumerate}
We set $W:=\Hom(G,G(n))$.
Let ${\frak G}(n):=Gr(V_m \otimes W,h[m](n))$
be the Grassmannian parametrizing $h[m](n)$-dimensional quotient
spaces of $V_m \otimes W$.
For a quotient
$V_m \otimes G \to E(m) \in {\frak Q}$,
let $F$ be the kernel.
Then for a sufficiently large $n$, we get that
\begin{enumerate}
\item
$\Hom(G,F(n)) \otimes G \to F(n)$ is surjective and
\item
$\Ext^i(G,F(n))=0$, $i>0$. 
\end{enumerate}
Hence we get a quotient vector space   
$V_m \otimes \Hom(G,G(n)) \to \Hom(G,E(m+n))$
of $V_m \otimes \Hom(G,G(n))$.
Thus we get a morphism ${\frak Q} \to {\frak G}(n)$.
As in \cite{Mum:1}, we can show that this morphism is a 
closed immersion
\begin{equation}
 {\frak Q} \hookrightarrow {\frak G}(n).
\end{equation}
Indeed, let ${\cal S} \subset V_m \otimes W \otimes {\cal O}_{{\frak G}(n)}$
be the universal subbundle.
We set 
\begin{equation}
 {\cal E}:=\coker({\cal S} \otimes G(-n) \to
 V_m \otimes W \otimes G(-n) \otimes {\cal O}_{{\frak G}(n)} 
 \overset{ev}{\to} V_m \otimes G \otimes {\cal O}_{{\frak G}(n)}). 
\end{equation}
We take a flattening stratification 
${\frak G}(n)=\coprod_i {\frak G}(n)_i$ of ${\cal E}$ \cite[sect. 8]{Mum:1}.
We may assume that each ${\frak G}(n)_i$ are connected.
Let ${\frak G}(n)_{h[m]}$ be the union of ${\frak G}(n)_i$ 
such that
the $G$-twisted Hilbert polynomial of ${\cal E}_x, x \in {\frak G}(n)_{h[m]}$
is $h[m]$. Then ${\frak G}(n)_{h[m]}$ is isomorphic to ${\frak Q}$.

$SL(V_m)$ acts on ${\frak G}(n)$. 
Let $L:={\cal O}_{{\frak G}(n)}(1)$ be the 
tautological line bundle on ${\frak G}(n)$.
Then $L$ has an $SL(V_m)$-linearization.
We consider the GIT semi-stability with respect to 
$L$.
\begin{prop}\label{prop:ss}
Let $\alpha:V_m \otimes W \to A$
be a quotient corresponding to a point of ${\frak G}(n)$.
Then it is GIT semi-stable with respect to $L$
if and only if 
\begin{equation}
 \dim V_m \dim \alpha(V' \otimes W)-
 \dim V' \dim \alpha(V_m \otimes W) \geq 0
\end{equation}
for all non-zero subspaces $V'$ of $V_m$.
\end{prop}

\begin{prop}
There is an integer $m_1$ such that for all $m \geq m_1$,
${\frak Q}^{ss}$ is contained in ${\frak G}(n)^{ss}$,
where $n \gg m$.
\end{prop}

\begin{proof}
We set 
\begin{equation}
{\cal F}:=\{E' \subset \widetilde{E}_q(m)|
E'=\im(V' \otimes G \to \widetilde{E}_q(m)), q \in {\frak Q}, 
V' \subset V_m \}.
\end{equation}
Since ${\cal F}$ is a bounded set, for a sufficiently large $n$ 
which depends on $m$,
\begin{enumerate}
\item
$\alpha(V' \otimes W)=\Hom(G,E'(n))$ and
\item
$\Ext^i(G,E'(n))=0$, $i>0$.
\end{enumerate}
Then $\dim \alpha(V' \otimes W)=\chi(G,E'(n))$.
Since $E$ is $G$-twisted semi-stable,
in the same way as in \cite[sect. 4]{Mar:4}, we see that
there is an integer $m_0$ such that for $m \geq m_0$ and
a subsheaf $E'$ of
$E(m)$,
\begin{equation}\label{eq:pgs}
 \frac{h^0(G,E')}{a_d(E')} \leq 
 \frac{h^0(G,E(m))}{a_d(E)}
\end{equation}
and the equality holds, if and only if 
\begin{equation}
 \frac{\chi(G,E'(n))}{a_d(E')}=\frac{\chi(G,E(m+n))}{a_d(E)}
\end{equation}
for all $n$.
Hence if the equality holds, then
$E'$ is $G$-twisted semi-stable and we may assume that ($\flat1,2$) holds
for $E'(-m)$.
In particular, $\dim V'=\chi(G,E')$.  
We set $m_1:=\max \{m_0,m(0) \}$.
For a sufficiently large $n \gg m$, we get
\begin{equation}\label{eq:a6}
 \left| \frac{\chi(G,E'(n))}{\chi(G,E(m+n))}
 -\frac{a_d(E')}{a_d(h)} \right|
 <\frac{1}{\dim V_m a_d(h)}.
\end{equation}
By \eqref{eq:a6}, 
if the inequality in \eqref{eq:pgs} is strict, then
we get 
\begin{equation}\label{eq:a7}
 \begin{split}
  &\dim V_m \dim \alpha(V' \otimes W)-
  \dim V' \dim \alpha(V_m \otimes W)\\
  > & \left(\dim V_m \frac{a_d(E')}{a_d(h)}-\dim V'-\frac{1}{a_d(h)}
  \right)\dim \alpha(V_m \otimes W)
  \geq 0.
 \end{split}
\end{equation}
If the equality holds in \eqref{eq:pgs}, then $\dim V'=\chi(G,E')$, 
and hence 
\begin{equation}
 \dim V_m \dim \alpha(V' \otimes W)-
  \dim V' \dim \alpha(V_m \otimes W)=0.
\end{equation}
Therefore our claim holds. 
\end{proof}

\begin{prop}\label{prop:proper}
There is an integer $m_2$ such that for all $m \geq m_2$,
${\frak Q}^{ss}$ is a closed subscheme of 
${\frak G}(n)^{ss}$,
where $n \gg m$.
\end{prop}

\begin{proof}
We choose an $m$ so that $h(m)/a_d(h)>1$.
We shall prove that 
${\frak Q}^{ss} \to {\frak G}(n)^{ss}$ is proper.
Let $(R,{\frak m})$ be a discrete valuation ring and $K$ the quotient field of $R$.
We set $T:=\Spec(R)$ and $U:=\Spec(K)$.
Let $U \to {\frak Q}^{ss}$ be a morphism such that $U \to {\frak Q}^{ss} \to 
{\frak G}(n)^{ss}$ is extended to
a morphism $T \to {\frak G}(n)^{ss}$.
Since ${\frak Q}$ is a closed subscheme of ${\frak G}(n)$,
there is a morphism $T \to {\frak Q}$, {\it i.e},
there is a flat family of quotients:
\begin{equation}
 V_m \otimes G \otimes {\cal O}_{T} \to 
 {\cal E}(m) \to 0.
\end{equation}
Let $\alpha:V_m \otimes W \otimes R \to 
\Hom_{p_{T*}}(G \otimes {\cal O}_{T},{\cal E}(m+n))$ 
be the quotient of
$V_m \otimes W \otimes R$ corresponding to
the morphism $T \to {\frak G}(n)^{ss}$.
We set $E:={\cal E} \otimes R/{\frak m}$.
\begin{claim}
$V_m \to \Hom(G,E(m))$ is injective.
\end{claim}
Indeed, we set $V':=\ker(V_m \to \Hom(G,E(m)))$.
Then $\alpha(V' \otimes W)=0$.
Hence we get

\begin{equation}
 \begin{split}
  0 \leq & \dim V_m \dim \alpha(V' \otimes W)-
  \dim V' \dim \alpha(V_m \otimes W)\\
  =& -\dim V' \dim \alpha(V_m \otimes W) \leq 0.
   \end{split}
\end{equation}
Therefore $V'=0$.

\begin{claim}\label{claim:2}
There is a rational number $\lambda$ which depends on $h$ 
such that $E$ is of type $\lambda$.
\end{claim}
Proof of the claim:
Let $E \to E''$ be a quotient of $E$.
Let $E'$ be the kernel of $E \to E''$.
We note that $V_m \to \Hom(G,E(m))$ is injective. 
We set $V':=V_m \cap \Hom(G,E'(m))$.
Then $h^0(G,E''(m)) \geq \dim V_m -\dim V'$.
Let $F$ be a subsheaf of $E(m)$ generated by $V'$.
Then $F$ belongs to ${\cal F}$.  
We set 
\begin{equation}
\varepsilon:=\frac{1}{a_d(h)! h(m)}.
\end{equation}
%
%$h(m)/a_d(h)- \varepsilon>0$.
%
Since ${\cal F}$ is a bounded set, for a sufficiently large $n$ 
which depends on $m$ and $\varepsilon$,
we have $\alpha(V' \otimes W)=\Hom(G,F(n))$,
$\Ext^i(G,F(n))=0$, $i>0$ and 
\begin{equation}\label{eq:a1}
 \left|\frac{ \dim \alpha(V' \otimes W)}{\dim \alpha(V_m \otimes W)}
  -\frac{a_d(F)}{a_d(h)} \right|<
 \varepsilon.
\end{equation}
Since $a_d(E') \geq a_d(F)$,
\begin{equation}
 \begin{split}
   \frac{h^0(G,E''(m))}{a_d(E'')} & \geq
   \frac{\dim V_m-\dim V'}{a_d(E'')}\\
  & > h(m)\left(\frac{a_d(h)-a_d(F)}{a_d(h)}\frac{1}{a_d(E'')}-
   \frac{\varepsilon}{a_d(E'')} \right)\\
  & > h(m)\left(\frac{a_d(h)-a_d(E')}{a_d(h)}\frac{1}{a_d(E'')}-
    \frac{\varepsilon}{a_d(E'')} \right)\\
  & \geq h(m)\left(\frac{1}{a_d(h)}-\varepsilon \right)>0.
 \end{split}
\end{equation}
There is a rational number $\lambda_1$ 
and an integer $m_3 \geq \lambda_1-a_{d-1}(h)/a_d(h)$ which depend on
$h(x)$ such that
\begin{equation}\label{eq:a5}
 h(m) \left(\frac{1}{a_d(h)}-\varepsilon \right) 
 \geq \frac{h(m)}{a_d(h)}-\frac{1}{a_d(h)!} \geq  
 \frac{1}{d!}\left(m+\frac{a_{d-1}(h)}{a_d(h)}-\lambda_1 \right)^d
\end{equation}
 for $m \geq m_3$.

By Lemma \ref{lem:valuative}, 
there is a purely $d$-dimensional sheaf $F$ with the
$G$-twisted Hilbert polynomial $h(x)$
and 
a map $E \to F$ whose kernel is a coherent sheaf of dimension less than $d$.
Let $F \to F''$ be a quotient such that $F''$ is $G$-twisted semi-stable.
We set $E':=\ker(E \to F'')$ and $E'':=\im(E \to F'')$.
Since $h^0(G,E''(m)) \leq h^0(G,F''(m))$ and $a_d(E'')=a_d(F'')$,
\begin{equation}
 \begin{split}
  \frac{h^0(G,F''(m))}{a_d(F'')} & \geq \frac{h^0(G,E''(m))}{a_d(E'')} \\    
  & \geq h(m) \left(\frac{1}{a_d(h)}-\varepsilon \right) >0.
 \end{split}
\end{equation}
Since $F''$ is $G$-twisted semi-stable, 
Lemma \ref{lem:section} implies that 
\begin{equation}
 \frac{h^0(G,F''(m))}{a_d(F'')} \leq 
   \left[\frac{1}{d!}\left(m+\frac{a_{d-1}(F'')}{a_d(F'')}+c \right)^d
   \right]_+
   \end{equation}
where $c$ is a constant which only depends on $a_d(h)$.
Since $h^0(G,F''(m))>0$, 
we get 
\begin{equation}
\frac{1}{d!}\left(m+\frac{a_{d-1}(F'')}{a_d(F'')}+c \right)^d
\geq 0
\end{equation}
 and
\begin{equation}
h(m) \left(\frac{1}{a_d(h)}-\varepsilon \right) \leq
 \frac{h^0(G,F''(m))}{a_d(F'')} \leq 
\frac{1}{d!}\left(m+\frac{a_{d-1}(F'')}{a_d(F'')}+c \right)^d.
\end{equation}
For $m \geq m_3$, \eqref{eq:a5} implies that
\begin{equation}
 \frac{a_{d-1}(h)}{a_d(h)}-\lambda_1 \leq 
 \frac{a_{d-1}(F'')}{a_d(F'')}+c.
\end{equation}
Hence
there is a rational number $\lambda$ which depends on
$h(x),\lambda_1$ such that
$F$ is of type $\lambda$.
Replacing $m$, we may assume that for all type $\lambda$ sheaves $I$
with the $G$-twisted Hilbert polynomial $h(x)$, 
\begin{enumerate}
\item
$\Hom(G,I(m)) \otimes G \to I(m)$ is surjective 
and
\item $\Ext^i(G,I(m))=0$, $i>0$.
\end{enumerate}
In particular $h^0(G,F(m))=h(m)=\dim V_m$.
Assume that $\Hom(G,E(m)) \to \Hom(G,F(m))$ is not injective and let
$V'$ be the kernel.
$J:=\im(V' \otimes G \to E(m))$ is of dimension less than $d$.
Hence we get $a_d(J)=0$.
By the inequality \eqref{eq:a1} and Proposition \ref{prop:ss},
we get a contradiction.
Thus $\Hom(G,E(m)) \to \Hom(G,F(m))$ is injective, and hence it is
isomorphic.
Since $\Hom(G,F(m)) \otimes G \to F(m)$ is surjective,
$E \to F$ must be surjective, which implies that it is isomorphic.
Therefore $E$ is of pure dimension $d$, of type $\lambda$ and 
$V_m \to \Hom(G,E(m))$ is an isomorphism.
Thus we complete the proof of Claim \ref{claim:2}. 

We set $m_2:=\max\{m_3,m(\lambda)\}$. 
Assume that there is a quotient $E \to E''$
which destabilizes the $G$-twisted semi-stability.
Since $E':=\ker (E \to E'')$
satisfies $(\flat 1, 2)$, we get that
$V'=\Hom(G,E'(m))$ and  
\begin{equation}
 \frac{\chi(G,E''(m))}{a_d(E'')} > 
 \frac{\chi(G,E(m))}{a_d(h)}-\varepsilon h(m).
\end{equation}
By the definition of $\varepsilon$,
we get  
\begin{equation}\label{eq:ss}
 \frac{\chi(G,E''(m))}{a_d(E'')} \geq 
 \frac{\chi(G,E(m))}{a_d(h)},
\end{equation}
which is a contradiction.
Therefore $E$ is $G$-twisted semi-stable.
Thus we get a lifting $T \to {\frak Q}^{ss}$ and 
conclude that ${\frak Q}^{ss} \to {\frak G}(n)^{ss}$ is proper.
\end{proof}
By standard arguments, we see that 
$SL(V_m)s, s \in {\frak Q}^{ss}$ is a closed orbit
if and only if 
the corresponding $G$-twisted semi-stable sheaf $E$
is isomorphic to
$\bigoplus_i E_i$,
where $E_i$ are $G$-twisted stable sheaves.

\end{NB}

\section{Twisted sheaves on a projective $K3$ surface}\label{sect:k3}

\subsection{Basic properties}
Let $X$ be a projective $K3$ surface and $p:Y \to X$ a projective bundle.

\begin{lem}\label{lem:equiv}
For a locally free $Y$-sheaf $E$,
$c_2({\bf R}p_*(E^{\vee} \otimes E)) \equiv
-(r-1)(w(E)^2) \mod 2r$.
\end{lem}

\begin{proof}
First we note that 
$(r-1)(D^2) \mod 2r$ is well-defined
for $D \in H^2(Z,\mu_r)$, $Z=X,Y$.
We take a representative $\alpha \in H^2(X,{\Bbb Z})$ of $w(E)$.
Then $c_1(E) \equiv p^*(\alpha) \mod r$.
Hence $c_2(p^*({\bf R}p_*(E^{\vee} \otimes E)))
=2rc_2(E)-(r-1)(c_1(E)^2) \equiv -(r-1)(p^*(\alpha^2)) \mod 2r$.
Since $H^4(X,{\Bbb Z})$ is a direct summand of
$H^4(Y,{\Bbb Z})$,
$c_2({\bf R}p_*(E^{\vee} \otimes E))
\equiv -(r-1)(\alpha^2) \mod 2r$.
\end{proof}

\begin{NB}
 
Let $D$ be a proper subscheme of $X$.
Since $\dim D \leq 1$,
$H^2(D,{\cal O}_D^{\times})=0$.
Hence there is a set of $\delta_{i,j} \in
H^0(U_i \cap U_j,{\cal O}_D^{\times})$ such that
${\cal O}_D^{\alpha}:=
\{({\cal O}_{D|U_i},\delta_{i,j})\}$ is an $\alpha$-twisted line bundle
on $D$.
\end{NB}

Let $K(X,Y)$ be the Grothendieck group of $Y$-sheaves.
\begin{lem}
\begin{enumerate}
\item[(1)]
There is a locally free $Y$-sheaf $E_0$
such that 
$\rk E_0=\min\{\rk E>0| E \in \Coh(X,Y)\}$.
\item[(2)]
$K(X,Y)={\Bbb Z}E_0 \oplus K(X,Y)_{\leq 1}$, where
$K(X,Y)_{\leq 1}$ is the submodule of $K(X,Y)$ generated by
$E \in \Coh(X,Y)$ of $\dim E \leq 1$.
\end{enumerate}
\end{lem}

\begin{proof}
(1) Let $F$ be a $Y$-sheaf such that
$\rk F=\min\{\rk E>0| E \in \Coh(X,Y)\}$.
Then $E_0:=F^{\vee \vee}$ satisfies the required properties.
(2)
We shall show that the image of $E \in \Coh(X,Y)$ in $K(X,Y)$
belongs to ${\Bbb Z}E_0 \oplus K(X,Y)_{\leq 1}$ by the induction of
$\rk E$. 
We may assume that $\rk E>0$.
Let $T$ be the torsion submodule of $E$. Then
$E=T+E/T$ in $K(X,Y)$.
 Since $\Hom(E_0(-nH),E/T) \ne 0$ for $n \gg 0$,
we have a non-zero homomorphism $\varphi:E_0(-nH) \to E/T$.
By our choice of $E_0$, $\varphi$ is injective.
Since $E_0(-nH)=E_0-E_{0|nH}$ in $K(X,Y)$,
$E=((E/T)/E_0+E_0)+(T-E_{0|nH})$. 
Since $\rk (E/T)/E_0<\rk E$,
we get $(E/T)/E_0 \in {\Bbb Z}E_0 \oplus K(X,Y)_{\leq 1}$, and hence
 $E$ also belongs to ${\Bbb Z}E_0 \oplus K(X,Y)_{\leq 1}$.
 \end{proof}

\begin{rem}
$\rk E_0$ is the order of the Brauer class of $Y$.
\end{rem}

\begin{NB}
$\chi(G,{\cal O}_D^{\alpha})=\chi(G^{\vee}_{|D}) 
\equiv (w(Y),D) \mod r$ and 
$\chi(G,k_P^{\alpha})=r$.
$(w(Y),D)=w(Y)_{|D} \equiv c_1({\bf R}p_*({\cal O}_D^{\alpha} \otimes
G^{\vee})) \mod r$.
\end{NB}

Let $\langle\;\;,\;\;\rangle$ be the Mukai pairing on $H^*(X,{\Bbb Z})$:
\begin{equation}
\langle x,y \rangle=-\int_X x^{\vee} y,\quad x,y \in H^*(X,{\Bbb Z}).
\end{equation} 
\begin{defn}
Let $G$ be a locally free $Y$-sheaf.
For a $Y$-sheaf $E$, we define a Mukai vector of $E$ as
\begin{equation}\label{eq:v-def}
\begin{split}
v_G(E):=& \frac{\ch({\bf R}p_*(E \otimes G^{\vee}))}
{\sqrt{\ch({\bf R}p_*(G \otimes G^{\vee}))}}\sqrt{\td_X} \\
=& (\rk(E),\zeta,b)
\in H^*(X,{\Bbb Q}),
\end{split}
\end{equation}
where $p^*(\zeta)=c_1(E)- \rk(E)\frac{c_1(G)}{\rk G}$ 
and $b \in {\Bbb Q}$.
More generally, for $G \in \Coh(X,Y)$ with $\rk G>0$,
we define $v_G(E)$ by \eqref{eq:v-def}.
\end{defn}

Since ${\bf R}p_*(E_1 \otimes G^{\vee}) \otimes 
{\bf R}p_*(E_2 \otimes G^{\vee})^{\vee} =
{\bf R}p_*(E_1 \otimes E_2^{\vee}) \otimes 
{\bf R}p_*(G \otimes G^{\vee})$,
\begin{equation}
\begin{split}
\langle v_G(E_1),v_G(E_2) \rangle=&
-\int_X \frac{\ch({\bf R}p_*(E_1 \otimes G^{\vee}))
\ch({\bf R}p_*(E_2 \otimes G^{\vee}))^{\vee}}
{\ch({\bf R}p_*(G \otimes G^{\vee}))}\td_X\\
=& -\int_X \ch({\bf R}p_*(E_1 \otimes E_2^{\vee}))\td_X\\
=& -\chi(E_2,E_1).
\end{split}
\end{equation}

\begin{NB}
\begin{equation}
\begin{split}
\frac{\sqrt{(r+\xi+b \omega)(r-\xi+b \omega)}}{(r+\xi+b \omega)}
&=\sqrt{\frac{(r-\xi+b \omega)}{(r+\xi+b \omega)}}\\
&=\sqrt{1-2\frac{\xi}{r}+2\left(\frac{\xi}{r}\right)^2}\\
&=e^{-\frac{\xi}{r}}
\end{split}
\end{equation}
\end{NB}

We define an integral structure on $H^*(X,{\Bbb Q})$ such that
$v_G(E)$ is integral. This is due to Huybrechts and Stellari
\cite{H-S:2}.
For a positive integer $r$ and $\xi \in H^2(X,{\Bbb Z})$,
we consider an injective homomorphism
\begin{equation}
\begin{matrix}
T_{-\xi/r}:& H^*(X,{\Bbb Z}) & \to & H^*(X,{\Bbb Q})\\
& x & \mapsto & e^{-\xi/r}x.
\end{matrix}
\end{equation}
$T_{-\xi/r}$ preserves the bilinear form $\langle\;\;,\;\; \rangle$.

\begin{lem}\label{lem:integral}
We take a representative $\xi \in H^2(X,{\Bbb Z})$ of
$w(G) \in H^2(X,\mu_{r})$, where $\rk(G)=r$.
We set $(\rk(E),D,a):=e^{\xi/r}v_G(E)$.
Then $(\rk(E),D,a)$ belongs to $H^*(X,{\Bbb Z})$
and $[D \mod \rk(E)]=w(E)$.
\end{lem}
\begin{proof}
We set $\sigma:=(c_1(G)-p^*(\xi))/r \in H^2(Y,{\Bbb Z})$.
Since $p^*(D)=p^*(\zeta)+\rk(E)p^*(\xi)/\rk(G)=
c_1(E)-\rk (E) \sigma \in H^2(Y,{\Bbb Z})$,
we get $D \in H^2(X,{\Bbb Z})$.
By Lemma \ref{lem:equiv}, we see that
\begin{equation}
\begin{split}
\langle e^{\xi/r}v_G(E),e^{\xi/r}v_G(E) \rangle=& 
\langle v_G(E),v_G(E) \rangle \\
= & c_2({\bf R}p_*(E \otimes E^{\vee}))-2\rk(E)^2 \\
\equiv &
(D^2) \mod 2 \rk (E).
\end{split}
\end{equation}
Hence $a \in {\Bbb Z}$. The last claim is obvious.
\end{proof}

\begin{rem}
$e^{\xi/r}v_G(E)$ is the same as the Mukai vector defined by
the rational $B$-field $\xi/r$ in \cite{H-S:2}.
More precisely,
there is a topological line bundle $L$ on $Y$ with
$c_1(L)=\sigma$ and $E \otimes L^{-1}$ is the pull-back of
a topological sheaf $E_{\xi/r}$ on $X$. Then we see that 
$e^{\xi/r}v_G(E)=\ch(E_{\xi/r})\sqrt{\td_X}$
(we use $H^i(X,{\Bbb Q})=0$ for $i>4$, or 
we deform $X$ so that $L$ becomes holomorphic).

\begin{NB}
Use $\ch(G \otimes L^{-1})/\ch((G \otimes L^{-1})^{-1})=e^{2\xi/r}$.
\end{NB} 
\end{rem}

\begin{defn}\cite{H-S:2}
We define a weight 2 Hodge structure on the lattice
$(H^*(X,{\Bbb Z}),\langle\;\;,\;\; \rangle)$
as 
\begin{equation}
\begin{split}
H^{2,0}(H^*(X,{\Bbb Z}) \otimes {\Bbb C}):= & 
T_{-\xi/r}^{-1}(H^{2,0}(X))\\ 
H^{1,1}(H^*(X,{\Bbb Z}) \otimes {\Bbb C}):=&
T_{-\xi/r}^{-1}(\bigoplus_{p=0}^2 H^{p,p}(X))\\
H^{0,2}(H^*(X,{\Bbb Z}) \otimes {\Bbb C}):=& 
T_{-\xi/r}^{-1}(H^{0,2}(X)).
\end{split}
\end{equation}
We denote this polarized Hodge structure by
$(H^*(X,{\Bbb Z}),\langle\;\;,\;\; \rangle,-\frac{\xi}{r})$.
\end{defn}

\begin{lem}
The Hodge structure
$(H^*(X,{\Bbb Z}),\langle\;\;,\;\; \rangle,-\frac{\xi}{r})$
depends only on the Brauer class $\delta'([\xi \mod r])$.
\end{lem}

\begin{proof}
If $\delta'([\xi \mod r])=\delta'([\xi' \mod r']) 
\in H^2(X,{\cal O}_X^{\times})$, then we have 
$r'\xi-r \xi'=L+rr' N$, where $L \in \NS(X)$ and
$N \in H^2(X,{\Bbb Z})$.
Then we have the following commutative diagram: 
\begin{equation}
\begin{CD}
H^*(X,{\Bbb Z}) @>{e^{-\frac{\xi}{r}}}>> H^*(X,{\Bbb Q})\\
@V{e^{-N}}VV @VV{e^{\frac{L}{rr'}}}V \\
H^*(X,{\Bbb Z}) @>>{e^{-\frac{\xi'}{r'}}}> H^*(X,{\Bbb Q}).
\end{CD}
\end{equation}
Thus we have an isometry of Hodge structures
\begin{equation}
(H^*(X,{\Bbb Z}),\langle\;\;,\;\; \rangle,-\frac{\xi}{r}) \cong
(H^*(X,{\Bbb Z}),\langle\;\;,\;\; \rangle,-\frac{\xi'}{r'}).
\end{equation}
\end{proof}

\begin{defn}
Let $Y \to X$ be a projective bundle and $G$ a locally free 
$Y$-sheaf. Let $\xi \in H^2(X,{\Bbb Z})$ be a lifting of
$w(G) \in H^2(X,\mu_r)$, where $r=\rk(G)$. 
\begin{enumerate}
\item
We define an integral Hodge structure of $H^*(X,{\Bbb Q})$
as $T_{-\xi/r}((H^*(X,{\Bbb Z}),\langle\;\;,\;\; \rangle,-\frac{\xi}{r}))$.
\item
$v:=(r,\zeta,b)$ is a Mukai vector, if 
$v \in 
T_{-\xi/r}(H^*(X,{\Bbb Z}))$
 and $\zeta \in \Pic(X) \otimes{\Bbb Q}$. 
 Moreover if $v$ is primitive in $T_{-\xi/r}(H^*(X,{\Bbb Z}))$,  
 then $v$ is primitive. 
\end{enumerate}
\end{defn}

\begin{defn}
Let $v:=(r,\zeta,b) \in H^*(X,{\Bbb Q})$ be a Mukai vector.
\begin{enumerate}
\item
$\overline{M}^{Y,G}_H(r,\zeta,b)$ (resp. ${M}^{Y,G}_H(r,\zeta,b)$)
denotes the coarse moduli space of $S$-equivalence classes of
$G$-twisted semi-stable (resp. stable) $Y$-sheaves $E$ with $v_G(E)=v$.
\item
${\cal M}^{Y,G}_H(r,\zeta,b)^{ss}$ (resp. ${\cal M}^{Y,G}_H(r,\zeta,b)^s$)
denotes the moduli stack of 
$G$-twisted semi-stable (resp. stable) $Y$-sheaves $E$ with $v_G(E)=v$.
\end{enumerate}
\end{defn}

\begin{lem}\label{lem:isom-G}
Assume that $o(w(Y))=o(w(Y'))$.
Then $\Xi_{Y \to Y'}^L$ induces an isomorphism
${\cal M}_H^{Y,G}(v)^{ss} \cong {\cal M}_H^{Y',G'}(v)^{ss}$,
where $G':=\Xi_{Y \to Y'}^L(G)$.
Moreover if $\dim Y=\dim Y'$ and $w(Y)=w(Y')$, then
${\cal M}_H^{Y,\epsilon(Y)}(v)^{ss} 
\cong {\cal M}_H^{Y',\epsilon(Y')}(v)^{ss}$.
\end{lem}
 
 \begin{proof}
 We use the notation in Lemma \ref{lem:o(w(Y))}.
 For a $Y$-sheaf $E$, we set $E':=\Xi_{Y \to Y'}^L(E)$.
 Then ${p'_Y}^*(E \otimes G^{\vee}) \cong p_{Y'}^*(E' \otimes {G'}^{\vee})$.
 Hence $v_G(E)=v_{G'}(E')$.
 If $\dim Y=\dim Y'$ and $w(Y)=w(Y')$, then 
since $w(\epsilon(Y))=w(\epsilon(Y'))$, replacing $L$ by
$L \otimes q^*(P)$, $P \in \Pic(X)$, we may assume that
$c_1(\Xi_{Y \to Y'}^L(\epsilon(Y)))=c_1(\epsilon(Y))$.
Thus $\Xi_{Y \to Y'}^L(\epsilon(Y))=\epsilon(Y)+T$ in 
$K(X,Y')$, where $T$ is a $Y$-sheaf with $\dim T=0$.   
From this fact, we get 
${\cal M}_H^{Y',\Xi_{Y \to Y'}^L(\epsilon(Y))}(v)^{ss}
={\cal M}_H^{Y',\epsilon(Y')}(v)^{ss}$.
 \end{proof}

\begin{NB}
If $\zeta \equiv \zeta' \mod r$ for $\zeta, \zeta' \in c_1(K(X,Y))$,
then
since $\zeta-\zeta' \in \NS(X)$,
$\zeta=\zeta'+rc_1(L)$, $L \in \Pic(X)$.
Hence if $H$ is a general polarization, then
$\overline{M}^{Y,G}_H(r,\zeta,b) \cong \overline{M}^{Y,G}_H(r,\zeta',b)$.  
\end{NB}

Let $E$ be a $Y$-sheaf. Then the Zariski tangent space 
of the Kuranishi space is $\Ext^1(E,E)$ and the obstruction space
is the kernel $\Ext^2(E,E)_0$
of the trace map
\begin{equation}
\tr:\Ext^2(E,E) \to H^2(Y,{\cal O}_Y) \cong H^2(X,{\cal O}_X).
\end{equation}
Hence as in the usual sheaves on a $K3$ surfaces
\cite{Mu:2}, we get the following.
\begin{prop}\label{prop:symplectic}
Let $E$ be a simple $Y$-sheaf. Then
the Kuranishi space is smooth of dimension $\langle v_G(E)^2 \rangle+2$
with a holomorphic symplectic form.
In particular, $\langle v_G(E)^2 \rangle \geq -2$.
\end{prop}

\begin{cor}\label{lem:bogomolov}
Let $E$ be a $\mu$-semi-stable $Y$-sheaf such that
$E=lE_0+F \in K(X,Y)$, $F \in K(X,Y)_{\leq 1}$.
Then $\langle v_G(E)^2 \rangle \geq -2 l^2$.
\end{cor}

\subsubsection{Wall and Chamber}
In this subsection, we generalize the notion of the wall
and the chamber for the usual stable sheaves
to the twisted case. 
\begin{lem}
Assume that there is an exact sequence of twisted sheaves
\begin{equation}\label{eq:seq2}
0 \to E_1 \to E \to E_2 \to 0,
\end{equation}
such that $E_i$, $i=1,2$ are $\mu$-semi-stable $Y$-sheaves.
We set $E_i=l_i E_0+F_i \in K(X,Y)$ with $F_i \in K(X,Y)_{\leq 1}$.
Then we have
\begin{equation}
\frac{\langle v_G(E)^2 \rangle}{l}+2l \geq 
-\frac{(l_2 v_G(F_1)-l_1 v_G(F_2))^2}{l l_1 l_2}.
\end{equation}
\end{lem}

This lemma easily follows from Corollary \ref{lem:bogomolov} 
and the following lemma.
\begin{lem}\label{lem:extension}
Let $E_0$ be a locally free $Y$-sheaf
such that 
$\rk E_0=\min\{\rk E>0| E \in \Coh(X,Y)\}$.
For an exact sequence of twisted sheaves
\begin{equation}\label{eq:seq}
0 \to E_1 \to E \to E_2 \to 0,
\end{equation}
we have
\begin{equation}
\frac{\langle v_G(E_1)^2 \rangle}{l_1}+
\frac{\langle v_G(E_2)^2 \rangle}{l_2}-
\frac{\langle v_G(E)^2 \rangle}{l}
=\frac{(l_2 v_G(F_1)-l_1 v_G(F_2))^2}{l l_1 l_2},
\end{equation}
where $E_i=l_i E_0+F_i$
and $E=l E_0+F$ in $K(X,Y)$ with $F_i,F \in K(X,Y)_{\leq 1}$.
%Then 
%$l_i:=\rk E_i/\rk E_0$, $i=1,2$ and $l=\rk E/\rk E_0$.
\end{lem}

\begin{proof}
%We set $E_i=l_i E_0+F_i$
%and $E=l E_0+F$ in $K(X,Y)$, where $F_i,F \in K(X,Y)_{\leq 1}$.
%Then 
\begin{equation}
\begin{split}
\frac{\langle v_G(E_1)^2 \rangle}{l_1}+
\frac{\langle v_G(E_2)^2 \rangle}{l_2}-
\frac{\langle v_G(E)^2 \rangle}{l} =&
\left(l_1 \langle v_G(E_0)^2 \rangle+2\langle v_G(E_0),v_G(F_1) \rangle
+\frac{\langle v_G(F_1),v_G(F_1) \rangle}{l_1}\right)\\
&+\left(l_2 \langle v_G(E_0)^2 \rangle+2\langle v_G(E_0),v_G(F_2) \rangle
+\frac{\langle v_G(F_2),v_G(F_2) \rangle}{l_2}\right)\\
& \quad-
\left(l \langle v_G(E_0)^2 \rangle+2\langle v_G(E_0),v_G(F) \rangle
+\frac{\langle v_G(F),v_G(F) \rangle}{l}\right)\\
=&
\frac{\langle v_G(F_1),v_G(F_1) \rangle}{l_1}+
\frac{\langle v_G(F_2),v_G(F_2) \rangle}{l_2}-
\frac{\langle v_G(F),v_G(F) \rangle}{l}\\
=&
\frac{(l_2 v_G(F_1)-l_1 v_G(F_2))^2}{l l_1 l_2}.
\end{split}
\end{equation}
\end{proof}

\begin{defn}
We set $v=v_G(lE_0+F)$, where $F$ is of dimension 1 or 0.
\begin{enumerate}
\item
For a
$\xi \in \NS(X)$ with 
$0<-(\xi^2) \leq \frac{l^2}{4}(2l^2+\langle v^2 \rangle)$,
we define a wall $W_{\xi}$ as
\begin{equation}
W_{\xi}:=\{L \in \Amp(X) \otimes {\Bbb R}|(\xi,L)=0 \}.
\end{equation}
\item
A chamber with respect to $v$ is a connected component of 
$\Amp(X) \otimes {\Bbb R} \setminus \bigcup_{\xi} W_{\xi}$.
\item
A polarization $H$ is general with respect to
$v$,
if $H$ does not lie on any wall.
\end{enumerate}
\end{defn}

\begin{rem}
The concept of chambers and walls are determined by $\rk(l E_0+F)$ and
$\langle v^2 \rangle$.
Thus they do not depend on the choice of $Y$ and $G$.
\end{rem} 

\begin{prop}\label{prop:chamber}
Keep notation as above.
\begin{enumerate}
\item
If $H$ and $H'$ belong to the same chamber, then
${\cal M}_H^{Y,G}(v)^{ss} \cong {\cal M}_{H'}^{Y,G}(v)^{ss}$.
\item
If $H$ is general, then 
${\cal M}_H^{Y,G}(v_G(F))^{ss} \cong 
{\cal M}_{H}^{Y,G'}(v_{G'}(F))^{ss}$
for $F \in K(X,Y)$ with $\rk F>0$.
\item
If 
\begin{equation}
\min\{-(D^2)>0| D \in \NS(X), (D,H)=0 \} >
 \frac{l^2}{4}(2l^2+\langle v^2 \rangle),
\end{equation}
then $H$ is general with respect to $v$.
\end{enumerate}
\end{prop}
The proof is standard (cf. \cite{H-L:1}) and is left to the reader.
%\begin{prop}\label{prop:chamber-general}
%If $H$ is general with respect to
%$v=(r,\zeta,b)$, then there is a projective bundle
%$Y' \to X$
% 
By Proposition \ref{prop:chamber} and
Proposition \ref{prop:symplectic}, we have 
%\begin{lem}
%Assume that $H$ is a general polarization with respect to $v$.
%For $G':=l'G+F'$, $l' \in {\Bbb Q}, F' \in K(X,Y)_{\leq 1} \otimes {\Bbb Q}$, 
%$E \in \Coh(X,Y)$ is $G$-twisted stable if and only if 
%$E$ is $G'$-twisted stable.  
%Thus 
%$M_H^{Y,G}(r,\xi,a) \cong M_H^{Y,G'}(r,0,b)$, where
%$G'$ is a locally free $Y$-sheaf of rank $r$ with 
%$c_1(G')=\xi+\frac{rc_1(G)}{\rk G}$ and $\xi^2-2ra=-2rb$.
%\end{lem}
%
%C\u ald\u araru \cite{C:2} developed a theory of derived category
%of twisted sheaves. 
%In particular, Grothendieck-Serre duality holds.
%Therefore the following follows by Mukai's arguments.
%
\begin{thm}\label{thm:symplectic}
Assume that $v$ is a primitive Mukai vector
and $H$ is general with respect to $v$.
Then all $G$-twisted semi-stable $Y$-sheaves $E$ with $v_G(E)=v$ are  
$G$-twisted stable.
In particular
$M_H^{Y,G}(v)$ is a projective manifold, if it is not empty.
\end{thm} 
In the next subsection, we show the non-emptyness of the moduli
space. We also show that $M_H^{Y,G}(v)$ is a $K3$ surface,
if $\langle v^2 \rangle=0$.

\begin{prop}(cf. \cite[Prop. 3.14]{Mu:4})\label{prop:simple}
Assume that $\Pic(X)={\Bbb Z}H$.
Let $E$ be a simple twisted sheaf with
$\langle v_G(E)^2 \rangle \leq 0$.
Then $E$ is stable.
\end{prop}

For the proof, we use Lemma \ref{lem:extension} and the following:

\begin{lem}\cite[Cor. 2.8]{Mu:4}
If $\Hom(E_1,E_2)=0$, then
\begin{equation}
\dim \Ext^1(E_1,E_1)+\dim \Ext^1(E_2,E_2) \leq
\dim \Ext^1(E,E).
\end{equation}
\end{lem}

\subsection{Existence of stable sheaves}

In this subsection, we shall show that the moduli space
of twisted sheaves is deformation equivalent
to the usual one.
In particular we show the non-emptyness of the moduli space.   

\begin{thm}\cite{H-S:1}\label{thm:brauer}
$H^1(X,PGL(r)) \to H^2(X,\mu_r)$ is surjective.
\end{thm}

\begin{prop}\label{prop:mu-stable}
For a $w \in H^2(X,\mu_r)$,
there is a ${\Bbb P}^{r-1}$-bundle
$p:Z \to X$ such that $w(Z)=w$ and $\epsilon(Z)$ is $\mu$-stable.
\end{prop}
D. Huybrechts informed us that the claim follows from the proof of
Theorem \ref{thm:brauer}.
Here we give another proof which works for other surfaces.
\begin{proof}
Let $p:Y \to X$ be a ${\Bbb P}^{r-1}$-bundle with
$w(Y)=w$.
We set $E_0:=\epsilon(Y)$.
In order to prove our claim,
it is sufficient to find a $\mu$-stable locally free $Y$-sheaf $E$
of rank $r$ with $c_1(E)=c_1(E_0)$.
For points $x_1,x_2,\dots,x_n \in X$,
let $F$ be a $Y$-sheaf which is the kernel of
a surjection $E_0 \to \bigoplus_{i=1}^n {\cal O}_{p^{-1}(x_i)}(1)$.
We take a smooth divisor $D \in |mH|$, $m \gg 0$. 
We set $\widetilde{D}:=p^{-1}(D)$.
Let $\Ext^i(F,F(-\widetilde{D}))_0$ be the kernel of the trace map
\begin{equation}
\Ext^i(F,F(-\widetilde{D})) \to H^i(Y,{\cal O}_Y(-\widetilde{D})) 
\cong H^i(X,{\cal O}_X(-D)).
\end{equation}
If $n \gg 0$, then the by the Serre duality,
 $\Ext^2(F,F(-\widetilde{D}))_0 \cong \Hom(F,F(\widetilde{D}))_0=0$.
\begin{NB}
Let $U$ be an affine scheme or an analytic open set and
$\phi:{\cal O}_U^{\oplus r} \to {\cal O}_U^{\oplus r}$ a homomorphism
such that $\tr{\phi}=0$. If $\phi(V)=V$ for all codimension 1
subspace $V \subset
{\cal O}_U^{\oplus r} \otimes k(x)$ and a general point 
$x \in U$, then $\phi=0$.
Hence if $\phi \ne 0$, then there is a quotient
$f:{\cal O}_U^{\oplus r} \to k(x)$ such that $\phi$ does not
preserve $\ker f$.     
\end{NB}
Hence $\Ext^1(F,F)_0 \to \Ext^1(F_{|\widetilde{D}},F_{|\widetilde{D}})_0$ 
is surjective.
Since $F_{|\widetilde{D}}$ deforms to a $\mu$-stable vector bundle
on $\widetilde{D}$,
$F$ deforms to a $Y$-sheaf $F'$ such that $F'_{|\widetilde{D}}$ is 
$\mu$-stable. Then $F'$ is also $\mu$-stable.
Then $E:=(F')^{\vee \vee}$ satisfies required properties.
\end{proof}

\begin{thm}\label{thm:deform}
Let $Y \to X$ be a projective bundle and
$G$ a locally free $Y$-sheaf.
Let $v_G:=(r,\zeta,b)$ be a primitive Mukai vector with $r>0$.
Then $M_H^{Y,G}(v_G)$ is an irreducible symplectic manifold
which is deformation equivalent to
$\Hilb_X^{\langle v_G^2 \rangle/2+1}$ for a general polarization $H$.
In particular 
\begin{enumerate}
\item[(1)]
$M_H^{Y,G}(v_G) \ne \emptyset$ if and only if
$\langle v_G^2 \rangle \geq -2$.
\item[(2)]
If $\langle v_G^2 \rangle=0$, then $M_H^{Y,G}(v_G)$ is a $K3$ surface.
\end{enumerate}
\end{thm}
We divide the proof into several steps.

Step 1 (Reduction to $M_H^{Y,\epsilon(Y)}(r,0,-a)$) :
Let $\xi$ be a lifting of $w(G)$.
Then $e^{\xi/\rk(G)}v_G=(r,D,b') \in H^*(X,{\Bbb Z})$.
By Theorem \ref{thm:brauer},
there is a projective bundle $Y' \to X$ such that
$w(Y')=[D \mod r]$.
Since $D/r-\xi/\rk(G)=\zeta/r \in \Pic(X) \otimes {\Bbb Q}$,
$o(w(Y'))=o(w(Y))$.
Let $G'$ be a locally free $Y$-sheaf such that
$\Xi_{Y \to Y'}^L(G')=\epsilon(Y')$, 
where we use the notation in Lemma \ref{lem:o(w(Y))}.
By Lemma \ref{lem:w(E)},
$w(G')=w(\epsilon(Y'))=[D \mod r]$.
Then replacing $L$ by $L \otimes q^*(P)$, $P \in \Pic(X)$,
we may assume that $e^{\xi/\rk G}v_G(G')=(r,D,c)$, $c \in {\Bbb Z}$.
Hence $v_{G'}(E)=(r,0,-a)$ for a $Y$-sheaf
$E$ with $v_G(E)=(r,\zeta,b)$. 
Since $H$ is general with respect to $(r,\zeta,b)$,
Proposition \ref{prop:chamber} implies that
$M_H^{Y,G}(r,\zeta,b) \cong M_H^{Y,G'}(r,0,-a)$.
By Lemma \ref{lem:isom-G},
$M_H^{Y,G'}(r,0,-a) \cong M_H^{Y',\epsilon(Y')}(r,0,-a)$.
Therefore replacing $(Y,G)$ by $(Y',\epsilon(Y'))$,
we shall prove the assertion for
$M_H^{Y,G}(r,0,-a)$ with $G=\epsilon(Y)$.

Step 2:
First we assume that $w(Y) \in \NS(X) \otimes \mu_r
\subset H^2(X,\mu_r)$.
Then the Brauer class of $Y$ is trivial, that is,
$Y={\Bbb P}(F)$ for a locally free sheaf $F$ on $X$.  
Since $H$ is general with respect to $(r,0,-a)$,
Proposition \ref{prop:chamber} (ii) and Lemma \ref{lem:isom-G} imply that
$M_H^{Y,G}(r,0,-a) \cong M_H^{X,{\cal O}_X}(r,D,c)$ 
with $2ra=(D^2)-2rc$.
By \cite[Thm. 8.1]{Y:7}, 
$M_H^{X,{\cal O}_X}(r,D,c)$ is deformation equivalent to
$\Hilb_X^{ra+1}$.

We next treat the general cases.
We shall deform the projective bundle $Y \to X$
to a projective bundle in Step 2.

Step 3: We first construct a local family of projective bundles.

\begin{prop}\label{prop:projective}
Let $f:({\cal X},{\cal H}) \to T$ be a family of polarized $K3$ surfaces.
Let $p:Y \to {\cal X}_{t_0}$ be a projective bundle
associated to a stable $Y$-sheaf $E$.
Then there is a smooth morphism $U \to T$ whose image contains $t_0$ 
and a projective bundle $p:{\cal Y} \to {\cal X} \times_T U$
such that ${\cal Y}_{t_0} \cong Y$.
\end{prop}

\begin{proof}
We note that $p_*(K_{Y/{\cal X}_{t_0}}^{\vee})$ is a vector bundle on 
${\cal X}_{t_0}$ and we have an embedding
$Y \hookrightarrow {\Bbb P}(p_*(K_{Y/{\cal X}_{t_0}}^{\vee}))$.
We take an embedding 
${\Bbb P}(p_*(K_{Y/{\cal X}_{t_0}}^{\vee})) \hookrightarrow
{\Bbb P}^{N-1} \times {\cal X}_{t_0}$
by a suitable quotient 
${\cal O}_{{\cal X}_{t_0}}(-n{\cal H}_{t_0})^{\oplus N} 
\to p_*(K_{Y/{\cal X}_{t_0}}^{\vee})$.
More generally, let ${\cal Y}_S \to {\cal X} \times_T S$
be a projective bundle and a surjective homomorphism
${\cal O}_{{\cal X} \times_T S}(-n{\cal H})^{\oplus N} 
\to p_*(K_{{\cal Y}_S/{\cal X} \times_T S}^{\vee})$.
Then we have an embedding
${\cal Y}_S \hookrightarrow {\Bbb P}^{N-1} \times {\cal X} \times_T S$.
 
Let ${\frak Y}$ be a connected component
of the Hilbert scheme $\Hilb_{{\Bbb P}^{N-1} \times {\cal X}/T}$
containing $Y$.
Let ${\cal Y} \subset {\Bbb P}^{N-1} \times {\cal X} \times_T {\frak Y}$
be the universal subscheme.
Let $\varphi:{\cal Y} \to {\cal X} \times_T {\frak Y}$
be the projection.
Let ${\frak Y}^0$ be an open subscheme of ${\frak Y}$
such that $\varphi_{|{\cal X} \times_T \{t\}}$ is smooth
and 
%$R^1 \varphi_*(T_{{\cal Y}/{\cal X} \times_T {\frak Y}})
%_{|{\cal X} \times_T \{t\}}=0$ for $t \in {\frak Y}^0$.
$H^1(T_{\varphi^{-1}(x,t)})=0$ for $(x,t) \in 
{\cal X} \times_T {\frak Y}^0$.
Since $Y \in {\frak Y}^0$, it is non-empty.
Then $\varphi$ is locally trivial on ${\cal X} \times_T {\frak Y}^0$.
Thus ${\cal Y} \to {\cal X} \times_T {\frak Y}^0$ is a projective bundle.  

If $Y$ is a projective bundle associated to a twisted vector bundle
$E$, then the obstruction for the infinitesimal liftings belongs to 
$H^2({\cal E}nd(E)/{\cal O}_X) \cong H^0({\cal E}nd(E)_0)^{\vee}$,
where ${\cal E}nd(E)_0$ is the trace free part of ${\cal E}nd(E)$.
\begin{NB}
By the Euler sequence, we have an exact sequence
$0 \to {\cal O}_X \to E^{\vee} \otimes E  \to p_*(T_{Y/X}) \to 0$.
Hence ${\bf R}p_*(T_{Y/X}) \cong {\cal E}nd(E)_0$.
\end{NB}
Hence if $E$ is simple (and $\rk E$ is not divisible by the characteristic), 
then there is no obstruction
for the infinitesimal liftings.
In particular ${\frak Y}^0 \to T$ is smooth at $Y$. 
\end{proof}

Step 4
{(A relative moduli space of twisted sheaves)}:
Let $f:({\cal X},{\cal H}) \to T$ be a family of polarized $K3$ surfaces
and $p:{\cal Y} \to {\cal X}$ a projective bundle on ${\cal X}$.
We set $g:=f \circ p$.
We note that 
$H^i({\cal Y}_t,\Omega_{{\cal Y}_t/{\cal X}_t})=0$, $i \ne 1$
and $H^1({\cal Y}_t,\Omega_{{\cal Y}_t/{\cal X}_t})={\Bbb C}$ 
for $t \in T$.
Hence  
$L:=\Ext^1_g(T_{{\cal Y}/{\cal X}},{\cal O}_{{\cal Y}}) \cong
R^1 g_*(\Omega_{{\cal Y}/{\cal X}})$ is a line bundle on $T$.
By the local-global spectral sequence, we have an isomorphism 
\begin{equation}
\Ext^1(T_{{\cal Y}/{\cal X}},g^*(L^{\vee})) \cong
H^0(T, \Ext^1_g(T_{{\cal Y}/{\cal X}},g^*(L^{\vee})))
\cong H^0(T,{\cal O}_T).
\end{equation}
We take the extension corresponding to $1 \in H^0(T,{\cal O}_T)$:
\begin{equation}
0 \to g^*(L^{\vee}) \to {\cal G} \to T_{{\cal Y}/{\cal X}} \to 0
\end{equation} 
such that ${\cal G}_t=\epsilon({\cal Y}_t)$.
Let $v:=(r,\zeta,b) \in R^* f_* {\Bbb Q}$ be a family of Mukai vectors
with $\zeta \in 
\NS({\cal X}/T) \otimes {\Bbb Q}$.
Then as in the absolute case,
we have a family of the moduli spaces of semi-stable twisted 
sheaves $\overline{M}_{ ({\cal X},{\cal H})/T}^{{\cal Y},\cal G}(v) \to T$
parametrizing ${\cal G}_t$-twisted semi-stable ${\cal Y}_t$-sheaves $E$ on
${\cal X}_t$, $t \in T$ with $v_{{\cal G}_t}(E)=v_t$.
$\overline{M}_{ ({\cal X},{\cal H})/T}^{{\cal Y},\cal G}(v) \to T$ 
is a projective morphism.
Let $E$ be a ${\cal G}_t$-twisted stable ${\cal Y}_t$-sheaf.
By our choice of $\zeta$, $\det(E)$ is unobstructed under deformations
over $T$, and hence $E$ itself is unobstructed.
Therefore ${M}_{ ({\cal X},{\cal H})/T}^{{\cal Y},\cal G}(v)$
is smooth over $T$.

Step 5 ({A family of $K3$ surfaces}):
Let 
${\cal M}_d$ be the moduli space of the 
polarized $K3$ surfaces $(X,H)$ with $(H^2)=2d$. 
${\cal M}_d$ is constructed as a quotient
of an open subscheme $T$ of a
suitable Hilbert scheme $\Hilb_{{\Bbb P}^N/{\Bbb C}}$.
Let $({\cal X},{\cal H}) \to T$ be the universal family.
Let $\Gamma$ be the abstruct $K3$ lattice and $h$ a primitive 
vector with $(h^2)=2d$.
Let ${\cal D}$ be the period domain for polarized $K3$ surfaces
$(X,H)$.
Let $\tau:\widetilde{T} \to T$ be the universal covering
and
$\phi_{\tilde{t}}:H^2({\cal X}_{\tau(\tilde{t})},{\Bbb Z}) 
\to \Gamma$, $\tilde{t} \in \widetilde{T}$ 
a trivialization on $\widetilde{T}$.
We may assume that $\phi_{\tilde{t}}({\cal H}_{\tau(\tilde{t})})=h$.
Then we have a period map ${\frak p}:\widetilde{T} \to {\cal D}$.
By the surjectivity of the period map, we can show that
${\frak p}$ is surjective:
Let $U$ be a suitable analytic neighborhood of a point $x \in {\cal D}$.
Then we have a family of polarized $K3$ surfaces
$({\cal X}_U,{\cal H}_U) \to U$ and an embedding of ${\cal X}$
as a subscheme of ${\Bbb P}^N \times U$.
Thus we have a morphism $h:U \to T$.
The embedding is unique up to the action
of $PGL(N+1)$.
Moreover if there is a point $\tilde{t}_0 \in \widetilde{T}$
such that ${\frak p}(\tilde{t}_0) \in U$, then
we have a lifting $\widetilde{h}:U \to \widetilde{T}$ of
$h:U \to T$ such that $\tilde{t}_0=\widetilde{h}({\frak p}(\tilde{t}_0))$.
Then
$U \to \widetilde{T} \to {\cal D}$
is the identity. 
Hence we can construct a lifting of any path on ${\cal D}$
intersecting ${\frak p}(\widetilde{T})$.
Since ${\cal D}$ is connected, we get the assertion.  
\begin{NB}
We may assume that $U$ is contractible.
If we take an embedding
${\cal X}_U \hookrightarrow {\Bbb P}^N \times U$, then
we have a diagram
\begin{equation}
\begin{CD}
({\cal X}_U,{\cal H}_U) @>>> {\Bbb P}^N \times U\\
@V{f}VV @|\\
({\cal X},{\cal H}) \times_{\frak H}U @>{h^*(i)}>> {\Bbb P}^N \times U\\
\end{CD}
\end{equation} 
Since we have an isomorphism
$\lambda:
({\cal X}_U,{\cal H}_U)_{{\frak p}(\tilde{t}_0)} \to 
({\cal X},{\cal H})_{\tau(\tilde{t}_0)}$,
we have another embedding
$({\cal X}_U)_{{\frak p}(\tilde{t}_0)} \to  {\cal X}_{\tau(\tilde{t}_0)}
\overset{i_{\tau(\tilde{t}_0)}}{\hookrightarrow} {\Bbb P}^N$.
Replacing $h$ by $g \circ h$, $g \in PGL(N+1)$,
we may assume that this map is the same as
$(h^*(i) \circ f)_{{\frak p}(\tilde{t}_0)}$.

*****
I don't need the following:
(We use the fact that an isomorphism $X \to X'$ preserving the
polarization can be lifted to the isomorphism of
the projective spaces.)
But use:
If $p_1,p_2:(X,L) \to {\Bbb P}^N$ be two embedding with
$p_1^*(H^0({\cal O}_{{\Bbb P}^N}(1))) \overset{\iota_1}{\to} H^0(L) \overset{\iota_2}{\leftarrow}
p_2^*(H^0({\cal O}_{{\Bbb P}^N}(1)))$, then
$\iota_2^{-1} \circ \iota_1$ induces an isomorphism 
$\iota:{\Bbb P}^N \to {\Bbb P}^N$ with
$\iota \circ p_1=p_2$. 
*****

Then we have
$h({\frak p}(\tilde{t}_0))=\tau(\tilde{t}_0)$ and
$\lambda=f_{{\frak p}(\tilde{t}_0)}$.
Let $\psi_u:H^2(({\cal X}_U)_u,{\Bbb Z}) \to \Gamma$ be the trivialization
on $U$. Then by the definition of ${\frak p}$,
$\psi_{{\frak p}(\tilde{t}_0)} \circ f_*^{-1}:
H^2( {\cal X}_{\tau(\tilde{t}_0)},{\Bbb Z}) \to 
H^2(({\cal X}_U)_{{\frak p}(\tilde{t}_0)},{\Bbb Z}) \to \Gamma$
coincides with $\phi_{\tilde{t}_0}$.
Hence if we take a local section $\sigma$ of $\tau$ in a neighborhood of
$\tau(\tilde{t}_0)$ with
$\sigma(\tau(\tilde{t}_0))=\tilde{t}_0$  
and set
$\tilde{h}:=\sigma \circ h$,
then $\phi_{\tilde{h}(u)}=\psi_u \circ f_*^{-1}$
in a neighborhood of ${\frak p}(\tilde{t}_0)$.
Since $U$ is contractible, $\tilde{h}$ is extended to
the whole of $U$ and we also have $\phi_{\tilde{h}(u)}=\psi_u \circ f_*^{-1}$.
Therefore ${\frak p} \circ \tilde{h}=1_U$.

*******
If we have an isomorphism $\xi:({\cal X}_U,{\cal H}_U)_u \to
({\cal X}_U,{\cal H}_U)_{u'}$ with a commutative diagram
\begin{equation}
\begin{CD}
H^2(({\cal X}_U)_u,{\Bbb Z}) @>{\psi_u}>> \Gamma \\
@V{\xi_*}VV @|\\
H^2(({\cal X}_U)_{u'},{\Bbb Z}) @>{\psi_{u'}}>> \Gamma, 
\end{CD}
\end{equation}
then $u'=u$, and hence $\xi_*=1$
(since the moduli space is fine).
*******
\end{NB}

\begin{NB}
Let $\gamma:[0,1] \to {\cal D}$ be a path in ${\cal D}$.    
Assume that there are local lifts $\widetilde{\gamma}_1$
on $[0,1/2]$ and $\widetilde{\gamma}_2$
on $[1/2,1]$ of $\gamma$.
Since $\widetilde{\gamma}_1(1/2)=g \circ \widetilde{\gamma}_2(1/2)$
for a $g \in PGL(N+1)$, replacing $\widetilde{\gamma}_2$ by
$g \circ \widetilde{\gamma}_2$, we have a lifting
of $\gamma$. 
For a general case, we cover a path $\gamma(t)$ by a finite number of
open covering such that there are local lifts on each open sets.
Then we can construct a lifting inductively.
\end{NB}

Step 6 (Reduction to step 2): 
We take a point $\widetilde{t} \in \widetilde{T}$.
We set
$(X,H):=
({\cal X}_{\tau(\widetilde{t})},{\cal H}_{\tau(\widetilde{t})})$. 
Let $p:Y \to X$ be a ${\Bbb P}^{r-1}$-bundle.
% such that
%$({\cal X}_{\tau(\widetilde{t})},{\cal H}_{\tau(\widetilde{t})})=(X,H)$.
Assume that $H$ is general with respect to $v:=(r,0,-a)$.
We take a $D \in \Gamma$ with $[D \mod r]=
\overline{\phi}_{\widetilde{t}}(w(Y))$.
Let $e_1,e_2,\dots,e_{22}$ be a ${\Bbb Z}$-basis of $\Gamma$ such that
$e_1=\phi_{\widetilde{t}}({\cal H}_{\tau(\widetilde{t})})$ and
$D=a e_1+b e_2$.
For an $\eta \in \bigoplus_{i=3}^{22} {\Bbb Z}e_i \subset
\Gamma$ with $(e_1^2)(\eta^2)-(e_1,\eta)^2<  0$,
we set
$\widetilde{\eta}:=e_2+rk \eta \in \Gamma$, $k \gg 0$. 
Since $\det
\begin{pmatrix}
(e_1^2) & (e_1,e_2+rk\eta)\\
(e_1,e_2+rk\eta) & ((e_2+rk\eta)^2)
\end{pmatrix}
 \ll 0$ for $k \gg 0$,
the signature of the primitive sublattice 
$L:={\Bbb Z}e_1 \oplus {\Bbb Z}\widetilde{\eta}$ of
$\Gamma$ is of type $(1,1)$.
Moreover $e_1^{\perp} \cap L$ 
does not contain a $(-2)$-vector.
We take a general $\omega \in L^{\perp} \cap \Gamma \otimes {\Bbb C}$ with
$(\omega,\omega)=0$ and $(\omega,\bar{\omega})>0$.
Then $\omega^{\perp} \cap \Gamma=L$. Replacing $\omega$ by its complex
conjugate if necessary, we may assume that
$\omega \in {\cal D}$.
Since ${\frak p}$ is surjective, there is a point
$\tilde{t}_1 \in \widetilde{{\frak H}}$ such that ${\frak p}(\tilde{t}_1)
=\omega$. 
Then ${\cal X}_{\tau(\widetilde{t}_1)}$ is a $K3$ surface with
$\Pic({\cal X}_{\tau(\widetilde{t}_1)})={\Bbb Z}
{\cal H}_{\tau(\widetilde{t}_1)} \oplus 
{\Bbb Z}\phi_{\widetilde{t}_1}^{-1}(e_2+rk\eta)$.
Hence $[\phi_{\widetilde{t}_1}^{-1}(D) \mod r]
=[\phi_{\widetilde{t}_1}^{-1}(a e_1+b\widetilde{\eta}) \mod r] \in 
\Pic({\cal X}_{\tau(\widetilde{t}_1)}) \otimes \mu_r$.
Since 
\begin{equation}
\min\{-(L^2)| 0 \ne L \in \Pic({\cal X}_{\tau(\widetilde{t}_1)}),
(L,{\cal H}_{\tau(\widetilde{t}_1)})=0 \} \gg 
\frac{r^2}{4}(2 r^2+\langle v^2 \rangle),
\end{equation}
Proposition \ref{prop:chamber} (iii) implies that
${\cal H}_{\tau(\widetilde{t}_1)}$ 
is a general polarization with respect to $v$. 
Then by the following lemma, we can reduce the proof to Step 2. 
Therefore we complete the proof of 
Theorem \ref{thm:deform}.

\begin{lem}\label{lem:deform}

For $\widetilde{t}_1,\widetilde{t}_2 \in \widetilde{T}$,
let $Y^i \to {\cal X}_{\tau(\tilde{t}_i)}$,
$i=1,2$ be ${\Bbb P}^{r-1}$-bundles
with $w(Y^i)=[\phi_{\tilde{t}_i}^{-1}(D) \mod r]$ and
$G_i:={\epsilon(Y^i)}$.
Let $v=(r,0,-a)$ be a primitive Mukai vector.
Assume that ${\cal H}_{\tau(\tilde{t}_i)}$, $i=1,2$ are general polarization.
Then
$M^{Y^1,G_1}_{{\cal H}_{\tau(\tilde{t}_1)}}(r,0,-a)$ is deformation
equivalent to
$M^{Y^2,G_2}_{{\cal H}_{\tau(\tilde{t}_2)}}(r,0,-a)$.
\end{lem}

\begin{proof} 
In order to simplify the notation, we denote
$M_{{\cal H}_t}^{Y,\epsilon(Y)}(r,0,-a)$ by
$M(Y)$ for a projective bundle $Y$ over $({\cal X}_t,{\cal H}_t)$.
By Proposition \ref{prop:mu-stable} and Lemma \ref{lem:isom-G},
we may assume that $\epsilon(Y^i)$ ($i=1,2$) is $\mu$-stable.
Let $\widetilde{\gamma}:[0,1] \to \widetilde{T}$ be a path
from $\tilde{t}_1=\widetilde{\gamma}(0)$ to 
$\tilde{t}_2=\widetilde{\gamma}(1)$ and
$\gamma:=\tau \circ \widetilde{\gamma}$.
Then we have a trivialization
$\overline{\phi}_s:H^2({\cal X}_{\gamma(s)},\mu_r) \to \Gamma \otimes_{\Bbb Z} \mu_r$.
By Proposition \ref{prop:mu-stable},
there is a projective bundle $Y_s \to {\cal X}_{\gamma(s)}$
such that $\overline{\phi}_s(w(Y_s))=[D \mod r]$
and $\epsilon(Y_s)$ is $\mu$-stable for each $s \in [0,1]$.
By Proposition \ref{prop:projective},
we have a family of projective bundles
${\cal Y}^s \to {\cal X} \times_T {\frak Y}^s$
over a $T$-scheme $\psi^s:{\frak Y}^s \to T$ 
such that there is a point 
$y^s \in {(\psi^s)}^{-1}(\gamma(s)) \subset {\frak Y}^s$ 
with $Y_s={\cal Y}^s_{y^s}$ and $\psi^s$ is smooth at $y^s$.
Then we have a family of moduli spaces
$\overline{M}_{ ({\cal X}\times_T {\frak Y}^s,\widetilde{\cal H})/{\frak Y}^s}
^{{\cal Y}^s,{\cal G}^s}(r,0,-a) \to  {\frak Y}^s$, where 
$\widetilde{\cal H}$ is the pull-back of ${\cal H}$
to ${\cal X}\times_T {\frak Y}^s$ (Step 4).
Since $\psi^s$ is smooth, $\psi^s({\frak Y}^s)$ is an open subscheme of
$T$ containing $\gamma(s)$.
We take an analytic open neighborhood $U_s$ of $\gamma(s)$ such that
$U_s$ is contractible and has a section $\sigma_s:U_s \to {\frak Y}^s$
with $\sigma_s(\gamma(s))=y^s$.
Let $V_s$ be a connected neighborhood of $s$
which is contained in $f^{-1}(U_s)$.
Since $[0,1]$ is compact, we can take a finite open covering of $[0,1]$:
$[0,1] = \cup_{j=1}^n V_{s_j}$, $s_1<s_2<\dots<s_n$.
Since $\{t \in T| \rk \Pic({\cal X}_t)=1 \}$ is a dense subset of $T$,
there is a point $t_j \in U_{s_j} \cap U_{s_{j+1}}$ such that
$t_j$ is sufficiently close to a point $\gamma(s_{j,j+1})$,
$s_{j,j+1} \in V_{s_j} \cap V_{s_{j+1}}$ and
$\Pic({\cal X}_{t_j})={\Bbb Z}{\cal H}_{t_j}$.
Under the identification
$H^2({\cal X}_t,\mu_r) \cong H^2({\cal X}_{\gamma(s)},\mu_r)$
for $t \in U_s$, we have
$w({\cal Y}^{s_j}_{\sigma_i(t_j)})=w({\cal Y}^{s_j}_{y^i})$
and $w({\cal Y}^{s_{j+1}}_{\sigma_{j+1}(t_j)})=
w({\cal Y}^{s_{j+1}}_{y^{j+1}})$.
Since $t_j$ is sufficiently close to the point $\gamma(s_{j,j+1})$,
%$s_{j,j+1} \in V_{s_j} \cap V_{s_{j+1}}$, 
we have
$w({\cal Y}^{s_j}_{\sigma_j(t_j)})=w({\cal Y}^{s_{j+1}}_{\sigma_{j+1}(t_j)})$.
Hence by Lemma \ref{lem:isom-G},
 $M({\cal Y}^{s_j}_{\sigma_j(t_j)})$
is isomorphic to
$M({{\cal Y}^{s_{j+1}}_{\sigma_{j+1}(t_j)}})$.
By Step 4,
$M({\cal Y}^{s_j}_{\sigma_j(t_{j-1})})$
is deformation equivalent to
$M({\cal Y}^{s_j}_{\sigma_j(t_{j})})$.
Therefore
$M({\cal Y}^{s_1}_{\sigma_1(t_{1})})$
is deformation equivalent to
$M({\cal Y}^{s_n}_{\sigma_n(t_{n-1})})$.
By using Step 4 again, we also see that
$M(Y^1)=M({\cal Y}^0_{y^0})$ is deformation equivalent to
$M({\cal Y}^{s_1}_{\sigma_1(t_{1})})$ and
$M(Y^2)=M({\cal Y}^1_{y^1})$ is deformation equivalent to
$M({\cal Y}^{s_n}_{\sigma_n(t_{n-1})})$.
Therefore our claim holds.
\end{proof}

\begin{NB}
The following proof is not correct (2004, Dec. 27).
\begin{proof} 
By Proposition \ref{prop:mu-stable} and Lemma \ref{lem:isom-G},
we may assume that $\epsilon(Y^i)$ ($i=1,2$) is $\mu$-stable.
Let
$\psi^i:{\frak Y}^i \to T$ be a morphism with a point 
$y^i \in {(\psi^i)}^{-1}({t}_i) \subset {\frak Y}^i$ and
${\cal Y}^i \to {\cal X} \times_T {\frak Y}^i$ 
a family of projective bundles such that
${\cal Y}^i_{y_0^i}=Y^i$ and $\psi^i$ is smooth 
(Proposition \ref{prop:projective}).
%Let $({\frak Y}^i \times_T \widetilde{T})_0$ be the connected component
%of ${\frak Y}^i \times_T \widetilde{T}$ containing 
%$(y_0^i,\tilde{t}_i)$.
%Since $\psi^i$ is smooth, $\psi^i({\frak Y}^i)$ is an open subscheme of
%$T$. Hence $\psi^i({\frak Y}^i) \times_T \widetilde{T}$ is connected.  
%Since ${\frak Y}^i \times_T \widetilde{T} \to \widetilde{T}$ is smooth,
%we see that $({\frak Y}^i \times_T \widetilde{T})_0 \to 
%\psi^i({\frak Y}^i) \times_T \widetilde{T}$ is surjective. 
We take a point $t \in \psi^1({\frak Y}^1) \cap
\psi^2({\frak Y}^2)$ with $\Pic({\cal X}_t)={\Bbb Z}{\cal H}_t$.
Let ${\cal G}^i$ be a vector bundle on ${\cal Y}^i$
with ${\cal G}^i_{y^i}=\epsilon({\cal Y}^i_{y^i}), y^i \in {\frak Y}^i$.
Let $\overline{\phi}_{\tilde{t}}:H^2({\cal X}_{\tau(\tilde{t})},\mu_r)
\to \Gamma \otimes_{\Bbb Z} \mu_r$ be the homomorphism
induced by $\phi_{\tilde{t}}$.
Since $w({\cal Y}^i_{y^i}) \in H^2({\cal X}_{\psi^i(y^i)},\mu_r)$
($y^i \in {\frak Y}^i$) is a section
of $R^2 f_{{\frak Y}^i *} \mu_r$, the map
$\gamma:
(y^i,\tilde{t}) \mapsto \overline{\phi}_{\tilde{t}}(w({\cal Y}^i_{y^i}))
\in \Gamma \otimes_{\Bbb Z} \mu_r$ is locally constant on
${\frak Y}^i \times_T \widetilde{T}$,
where $f_{{\frak Y}^i}$ is the base change of $f$
by $\psi^i$. 
By the definition of ${\frak Y}^i$,
$w({\cal Y}^i_{y^i})$ does not depend on the choice of $y^i \in {\frak Y}^i$.
Hence $\gamma$ is regarded as a locally constant function on
a connected manifold $\psi^i({\frak Y}^i) \times_T \widetilde{T}$. 
Therefore $\gamma$ is constant.
Thus
$\overline{\phi}_{\tilde{t}}(w({\cal Y}^i_{y^i}))=
\overline{\phi}_{\tilde{t}_i}(w({\cal Y}^i_{y_0^i}))=
[D \mod r]$ for $(y^i,\tilde{t}) \in {\frak Y}^i \times_T \widetilde{T}$.
In particular $w({\cal Y}^1_{y^1})=w({\cal Y}^2_{y^2}) 
\in H^2({\cal X}_t,\mu_r)$
for $y^i \in (\psi^i)^{-1}(t)$.
By Lemma \ref{lem:isom-G},
$M^{{\cal Y}^1_{y^1},{\cal G}_{y^1}}_{{\cal H}_{\psi^1(y^1)}}(r,0,-a) \cong
M^{{\cal Y}^2_{y^2},{\cal G}_{y^2}}_{{\cal H}_{\psi^2(y^2)}}(r,0,-a)$.
By Step 4, $M^{Y^i,G_i}_{{\cal H}_{t_i}}(r,0,-a)=
M^{{\cal Y}^i_{y_0^i},{\cal G}_{y_0^i}}_{{\cal H}_{\psi^i(y_0^i)}}(r,0,-a)$ 
is deformation equivalent to
$M^{{\cal Y}^i_{y^i},{\cal G}_{y^i}}_{{\cal H}_{\psi^i(y^i)}}(r,0,-a)$.
Hence we get our claim.
\end{proof}
\end{NB}

\begin{NB}
$\tau:\widetilde{T} \to T$ be the universal covering of $T$.
Then we have a trivialization 
$\phi_{\tilde{t}}:H^2({\cal X}_{\tau(\tilde{t})},{\Bbb Z})
\to \Gamma$, $\tilde{t} \in \widetilde{T}$.
We have a period map ${\frak p}:\widetilde{T} \to {\cal D}$.
$w({\cal Y})_h \in H^2({\cal X}_{\zeta(h)},\mu_r)$.
$(h,\tilde{t}) \in {\frak H} \times_T \widetilde{T}$,
$w({\cal Y})_h =[\phi_{\tilde{t}}^{-1}(D) \mod r]$.

For $t_1, t_2 \in \zeta({\frak H})$,
$M_{{\cal H}_{t_1}}^{G_1}(r,\xi_1,a)$ is deformation equivalent to
$M_{{\cal H}_{t_2}}^{G_2}(r,\xi_2,a)$, where $G_i$, $i=1.2$ 
are twisted locally free
sheaves such that 
$w({\Bbb P}(G_i^{\vee}))=[\phi_{\tilde{t_i}}^{-1}(D) \mod r]$
and $[\xi_i \mod r]=[\phi_{\tilde{t_i}}^{-1}(D) \mod r]$.
\end{NB}

\begin{rem}
Let $v_G:=(r,\zeta,b)$ be a Mukai vector with $r,\langle v_G^2 \rangle>0$
which is not necessary primitive.
By the same proof, we can also show that 
$\overline{M}_H^{Y,G}(v_G)$ is an irreducible normal variety
for a general $H$ (cf. \cite{Y:9}). 
\end{rem}

\subsection{The second cohomology groups of moduli spaces}

By Theorem \ref{thm:deform},
$M_H^{Y,G}(v_G)$ is an irreducible symplectic manifold,
if $v_G$ is primitive and $H$ is general.
Then $H^2(M_H^{Y,G}(v_G),{\Bbb Z})$ is equipped with a bilinear
form called the Beauville form.
In this subsection, we shall describe the Beauville
form in terms of the Mukai lattice.
   
Let $p:Y \to X$ be a projective bundle with $w(Y)=[\xi \mod r]$
and set $G:=\epsilon(Y)$.
We consider a Mukai lattice with a Hodge structure
$(H^*(X,{\Bbb Z}),\langle\;\;,\;\; \rangle,-\frac{\xi}{r})$
in this subsection.
We set $w:=r(1,0,\frac{a}{r}-\frac{1}{2}\frac{(\xi^2)}{r^2})$,
$a \in {\Bbb Z}$.
In this subsection, we assume that $w$ is primitive, that is,
$\gcd(r,\xi,a)=1$. 
We set $v:=w e^{\xi/r}=(r,\xi,a) \in H^*(X,{\Bbb Z})$.
Then $v$ is algebraic.

Let $q:\widetilde{M_H^{Y,G}(w)} \to M_H^{Y,G}(w)$ be a projective bundle
in subsection \ref{subsubsect:family}
and ${\cal E}$ the family of twisted sheaves on
$Y \times \widetilde{M_H^{Y,G}(w)}$.
We set $W^{\vee}:=\epsilon(\widetilde{M_H^{Y,G}(w)})$. 
Let $\widetilde{\pi}_{\widetilde{M_H^{Y,G}(w)}}:
Y \times \widetilde{M_H^{Y,G}(w)} \to \widetilde{M_H^{Y,G}(w)}$
and $\widetilde{\pi}_{Y}:
Y \times \widetilde{M_H^{Y,G}(w)} \to Y$ be projections.
Then
$(1_Y \times q)_*({\cal E} \otimes 
\widetilde{\pi}_{\widetilde{M_H^{Y,G}(w)}}^*(W^{\vee}))$ 
is a quasi-universal family on
$Y \times M_H^{Y,G}(w)$.

Let $\pi_X:X \times M_H^{Y,G}(w) \to X$ be the projection.
We define a homomorphism
$\theta_v^G:v^{\perp} \to H^*(M_H^{Y,G}(w),{\Bbb Q})$
by
\begin{equation}
\theta_v^G(u):=\int_{X} [{\cal Q}^{\vee} \pi_X^*(e^{-\xi/r}u)]_3
\end{equation}
where $[...]_3$ means the degree 6 part and
\begin{equation}
\begin{split}
{\cal Q}:=& 
\frac{\sqrt{\td_X}}{\sqrt{\ch({\bf R}p_*(G^{\vee} \otimes G))}}
\frac{\sqrt{\td_{M_H^{Y,G}(w)}}}{\sqrt{\ch({\bf R}q_*({W}^{\vee} \otimes W))}} 
\ch({\bf R}(p \times q)_*(\widetilde{\pi}_Y^*( G^{\vee}) 
\otimes {\cal E} \otimes 
\widetilde{\pi}_{\widetilde{M_H^{Y,G}(w)}}^*({W}^{\vee})))\\
& \in H^*(X \times M_H^{Y,G}(w),{\Bbb Q}).
\end{split}
\end{equation}

\begin{rem}\label{rem:integral}
If $\xi$ is algebraic, then $Y$ is isomorphic to the projective bundle
${\Bbb P}(F^{\vee})$ and $G=F^{\vee} \otimes {\cal O}_Y(1)$, 
where $F$ is a vector bundle of rank $r$ on $X$ with $c_1(F)=-\xi$.
In this case, $M_H^{Y,G}(w)$ is the usual moduli space of stable sheaves
$F$ with the Mukai vector $v$ and  
${\bf R}(p \times q)_*( \widetilde{\pi}_Y^*({\cal O}_Y(-1))
 \otimes {\cal E} \otimes
\widetilde{\pi}_{\widetilde{M_H^{Y,G}(w)}}^*({W}^{\vee}))$ 
is a quasi-universal family.
Since $\ch F/\sqrt{\ch(F \otimes F^{\vee})}=e^{-\xi/r}$, we have
\begin{equation}
{\cal Q}=
e^{-\frac{\xi}{r}}\sqrt{\td_X}\frac{\sqrt{\td_{M_H^{Y,G}(w)}}}{\sqrt{\ch({\bf R}q_*({W}^{\vee} \otimes W))}}
\ch({\bf R}(p \times q)_*(\widetilde{\pi}_Y^*( {\cal O}_Y(-1))
 \otimes {\cal E} \otimes 
\widetilde{\pi}_{\widetilde{M_H^{Y,G}(w)}}^*({W}^{\vee}))).
\end{equation}
Hence $\theta_v^G$ is the usual Mukai homomorphism,
which is defined over ${\Bbb Z}$.
\end{rem}

Let $p':Y' \to X$ be another ${\Bbb P}^{r-1}$-bundle with $w(Y')=w(Y)$.
Then by the proof of Lemma \ref{lem:isom-G},
we see that the following diagram is commutative:
\begin{equation}
\begin{CD}
v^{\perp}@= v^{\perp}\\
@V{\theta_v^G}VV @VV{\theta_v^{G'}}V\\
H^2(M_H^{Y,G}(w),{\Bbb Q}) @>>> H^2(M_H^{Y',G'}(w),{\Bbb Q}),
\end{CD}
\end{equation}
where $G':=\Xi_{Y \to Y'}^L(G)=\epsilon(Y')$.
Since ${\cal Q}$ is algebraic, 
$\theta_v^G$ preserves the Hodge structure.
By the deformation argument, 
Remark \ref{rem:integral} implies that
$\theta_v^G$ is defined over ${\Bbb Z}$.
Moreover it preserves the bilinear forms.

\begin{thm}
For $\xi \in H^2(X,{\Bbb Z})$ with $[\xi \mod r]=w(Y)$,
we set $v=w e^{\xi/r}$.
\begin{enumerate}
\item
If $\langle v^2 \rangle>0$, then
$\theta_v^G:v^{\perp} \to H^2(M_H^{Y,G}(w),{\Bbb Z})$ is an isometry
of the Hodge structures.
\item
If $\langle v^2 \rangle=0$, then $\theta_v^G$ induces an isometry
of the Hodge structures
$v^{\perp}/{\Bbb Z}v \to H^2(M_H^{Y,G}(w),{\Bbb Z})$.
\end{enumerate} 
\end{thm}
The second claim is due to Mukai \cite{Mu:5}.

\section{Fourier-Mukai transform}\label{sect:FM}
\begin{NB} 
Let $\mathrm{Mod}(X,Y)$ be the subcategory of 
the category of ${\cal O}_Y$-modules
consisting of $E$ with $p^*(p_*(E \otimes G^{\vee}))
\cong E \otimes G^{\vee}$.
Let $I$ be an injective sheaf on $X$.
Then $p^*(I) \otimes G$ is injective in $\mathrm{Mod}(X,Y)$.
Indeed $\Hom(E,p^*(I) \otimes G)\cong \Hom(E \otimes G^{\vee},p^*(I))
\cong \Hom(p_*(E \otimes G^{\vee}),I)$.
Hence $E \mapsto \Hom(E,p^*(I) \otimes G)$ is an exact functor.
Since there is an injection $E \to (E \otimes G^{\vee}) \otimes G$,
$\mathrm{Mod}(X,Y)$ is enough injective. 
\end{NB}

\subsection{Integral functor}

Let $p:Y \to X$ be a projective bundle such that
$\delta([Y])=[\alpha] \in \Br(X)$ and
 $p':Y' \to X'$ a projective bundle such that
$\delta([Y'])=[\alpha'] \in \Br(X')$.
Let ${\pi}_X:X' \times X \to X$
and ${\pi}_{X'}:X' \times X \to X'$ be projections.
We also let $\widetilde{\pi}_Y:Y' \times Y \to Y$
and $\widetilde{\pi}_{Y'}:Y' \times Y \to Y'$ be projections.
We set $G:=\epsilon(Y)$ and $G':=\epsilon(Y')$.
\begin{defn}
Let $\Coh(X' \times X,Y',Y)$ be the subcategory
of $\Coh(Y' \times Y)$ such that $Q \in \Coh(Y' \times Y)$
belongs to $\Coh(X' \times X,Y',Y)$ if and only if
$(p' \times p)^* (p' \times p)_*
(G' \otimes Q \otimes G^{\vee}) \cong 
G' \otimes Q \otimes G^{\vee}$.
In terms of local trivialization of $p,p'$, this is equivalent to
$Q_{|Y_i' \times Y_j} \cong
{\cal O}_{Y_i'}(-\lambda_i') \boxtimes {\cal O}_{Y_j}(\lambda_j)
\otimes (p' \times p)^*({Q}_{ij})$,
${Q}_{ij} \in \Coh(U_i' \times U_j)$.
$\Coh(X' \times X,Y',Y)$ is equivalent to
$\Coh(X' \times X,{\alpha'}^{-1} \times \alpha)$.
\end{defn}

\begin{rem}
We take twisted line bundles ${\cal L}({p'}^*({\alpha'}^{-1}))$ 
on $Y'$ and
${\cal L}(p^*(\alpha^{-1}))$ on $Y$ respectively 
which give equivalences
$\Lambda^{{\cal L}({p'}^*({\alpha'}^{-1}))}:
\Coh(X',Y') \cong \Coh(X',\alpha')$ and
$\Lambda^{{\cal L}(p^*({\alpha}^{-1}))}:
\Coh(X,Y) \cong \Coh(X,\alpha)$ in \eqref{eq:Y=X}.
Then we have an equivalence
$\Lambda^{{\cal L}({p'}^*({\alpha'}^{-1}))^{\vee}} \times 
\Lambda^{{\cal L}(p^*({\alpha}^{-1}))}$:
\begin{equation}
\begin{matrix}
\Coh(X' \times X,Y',Y)& \to &\Coh(X' \times X,{\alpha'}^{-1} \times \alpha)\\
 Q & \mapsto & 
(p' \times p)_*({\cal L}({p'}^*({\alpha'}^{-1})) \otimes
Q \otimes {\cal L}(p^*({\alpha}^{-1}))^{\vee}).
\end{matrix}
\end{equation}
\end{rem}

Let ${\bf D}(X' \times X,Y',Y) \cong 
{\bf D}(X' \times X,{\alpha'}^{-1} \times \alpha)$ be the bounded 
derived category of $\Coh(X' \times X,Y',Y)$.
\begin{NB}
${\cal Q} \in {\bf D}(X' \times X,{\alpha'}^{-1} \times \alpha)$ corresponds
an object $\widetilde{\cal Q} \in {\bf D}(Y' \times Y)$ 
which is represented by a complex
$\cdots \to \widetilde{\cal Q}^k \to \widetilde{\cal Q}^{k+1} \to \cdots$
with $(p' \times p)^* (p' \times p)_*
(G' \otimes \widetilde{\cal Q}^k \otimes G^{\vee}) \cong 
G' \otimes \widetilde{\cal Q}^k \otimes G^{\vee}$.
or ${\cal Q}_{|Y_i' \times Y_j}^k =
{\cal O}_{Y_i'}(-\lambda_i') \boxtimes {\cal O}_{Y_j}(\lambda_i)
\otimes (p' \times p)^*({\cal Q}^k_{ij})$,
${\cal Q}^k_{ij} \in \Coh(U_i' \times U_j)$.
\end{NB}
For ${\cal Q} \in {\bf D}(X' \times X,Y',Y)$, we define an integral
functor 
\begin{equation}
\begin{matrix}
\Phi_{X' \to X}^{\widetilde{\cal Q}}: 
&{\bf D}(X',Y')& \to &{\bf D}(X,Y)\\
& x & \mapsto & {\bf R}\widetilde{\pi}_{Y*}
({\cal Q} \otimes \widetilde{\pi}_{Y'}^*(x)).
\end{matrix}
\end{equation}

For ${\cal Q} \in {\bf D}(X' \times X,Y',Y)$ and
${\cal R} \in {\bf D}(X'' \times X', Y'',Y')$,
we have
\begin{equation}
\Phi_{X' \to X}^{{\cal Q}} \circ 
\Phi_{X'' \to X'}^{{\cal R}}=
\Phi_{X'' \to X}^{{\cal S}},
\end{equation}
where 
${\cal S}=
{\bf R}\widetilde{\pi}_{Y'' \times Y*}
(\widetilde{\pi}^*_{Y'' \times Y'}({\cal R}) \otimes
 \widetilde{\pi}^*_{Y' \times Y}({\cal Q}))$ and
$\widetilde{\pi}^*_{(\;\;)}:Y'' \times Y' \times Y \to (\;\;)$
is the projection.

\begin{NB}
For a torsion class $\alpha:=\{\alpha_{ijk}\}$,
let $p:Y \to X$ be a projective bundle associated to
$\alpha$.

${\cal Q} \in {\bf D}(X' \times X,\beta^{-1} \times \alpha)$ corresponds
an object ${\cal Q} \in {\bf D}(Y' \times Y)$ such that
${\cal Q}_{|Y_i' \times Y_j} =
{\cal O}_{Y_i'}(-\lambda_i') \boxtimes {\cal O}_{Y_j}(-\lambda_i)
\otimes (p' \times p)^*({\cal Q}_{ij})$.
\begin{equation}
\begin{matrix}
\Phi_{X' \to X}^{\cal Q}: &{\bf D}(X',\beta^{-1})& \to &{\bf D}(X,\alpha)\\
& x & \mapsto & {\bf R}\pi_{X*}({\cal E} \otimes \pi_{X'}^*(x)).
\end{matrix}
\end{equation}

For ${\cal Q} \in {\bf D}(X' \times X,\beta^{-1} \times \alpha)$ and
${\cal R} \in {\bf D}(X'' \times X', \gamma^{-1} \times \beta)$,
we have
\begin{equation}
\Phi_{X' \to X}^{\cal Q} \circ \Phi_{X'' \to X'}^{\cal R}=
\Phi_{X'' \to X}^{\cal S},
\end{equation}
where 
${\cal S}=
{\bf R}\pi_{Y'' \times Y*}({\cal R} \boxtimes_{{\cal O}_{Y'}}
 {\cal Q})$.

\end{NB}

\subsubsection{Cohomological correspondence}
For simplicity,
we denote the pull-backs of $G$ and $G'$ to $Y' \times Y$
by the same letters.
For example
$G' \otimes {\cal Q} \otimes G^{\vee}$ implies 
$\pi_{Y'}^*(G') \otimes {\cal Q} \otimes \pi_{Y}(G^{\vee})$.
\begin{NB}
We note that
\begin{equation}
\frac{1}{\sqrt{\ch({\bf R}p_*(G^{\vee} \otimes G))}} \in H^*(X,{\Bbb Q}).
\end{equation}

We set 
\begin{equation}
G' \boxtimes_{{\cal O}_{Y'}} {\cal Q} \boxtimes_{{\cal O}_{Y'}} G^{\vee}:=
\pi_{Y'}^*(G') \otimes {\cal Q} \otimes
\pi_{Y'}^*(G^{\vee}).
\end{equation}
\end{NB}
We note that
\begin{equation}
{\bf R}(p' \times p)_*(G' \otimes {\cal Q} \otimes G^{\vee}) \in 
{\bf D}(X' \times X)
\end{equation}  
satisfies
\begin{equation}
(p' \times p)^*({\bf R}(p' \times p)_*(G' \otimes {\cal Q} \otimes G^{\vee}))
=G' \otimes {\cal Q} \otimes G^{\vee}.
\end{equation}
We define a homomorphism
\begin{equation}
\Psi_{X' \to X}^{\cal Q}:
H^*(X',{\Bbb Q}) \to H^*(X,{\Bbb Q})
\end{equation}
by
\begin{equation}
\begin{split}
&\Psi_{X' \to X}^{\cal Q}(y)\\
:=&\pi_{X*} \circ (p' \times p)_*
\left((p' \times p)^* \circ \pi_{X'}^*(y)
 \ch({G'})\ch({\cal Q})\ch(G^{\vee})
\frac{\sqrt{\td_{X'}}\td_{Y'/X'}}{\sqrt{\ch({G'}^{\vee} \otimes G')}} 
\frac{\sqrt{\td_X}\td_{Y/X}}{\sqrt{\ch(G^{\vee} \otimes G)}}\right)\\
=&
\pi_{X*}\left( \pi_{X'}^*(y) 
\frac{\sqrt{\td_{X'}}}{\sqrt{\ch({\bf R}p'_*({G'}^{\vee} \otimes G'))}} 
\frac{\sqrt{\td_X}}{\sqrt{\ch({\bf R}p_*(G^{\vee} \otimes G))}}
\ch({\bf R}(p' \times p)_*({G'}  \otimes {\cal Q}\otimes G^{\vee}))\right),
\end{split}
\end{equation}
where $\td_X$, $\td_{X'}$,... are identified with their pull-backs.   

\begin{lem}
$\Psi_{X'' \to X}^{\cal S}=
\Psi_{X' \to X}^{\cal Q} \circ \Psi_{X'' \to X'}^{\cal R}$.
\end{lem}

\begin{proof}
$\pi_{(\;\;)}:X'' \times X' \times X \to (\;\;)$ 
denotes the projection to $(\;\;)$.
% and
%$\widetilde{\varpi}_{(\;\;)}:Y'' \times Y' \times Y \to (\;\;)$ 
%$the projection to $(\;\;)$.
We note that
\begin{equation}
\begin{split}
& \pi_{X'' \times X}^* \left({\bf R}(p'' \times p')_*
(G'' \otimes {\cal R} \otimes {G'}^{\vee})\right) \otimes
\pi_{X' \times X}^* \left({\bf R}(p' \times p)_*
(G' \otimes {\cal Q} \otimes {G}^{\vee})\right) \\
=& {\bf R}(p'' \times p' \times p)_*
(G'' \otimes {\cal R} \otimes {\cal Q} \otimes G^{\vee})
\otimes \pi_{X'}^*({\bf R}p'_*({G'}^{\vee} \otimes G')).
\end{split}
\end{equation}
Then
\begin{equation}
\begin{split}
& \pi_{X'' \times X}^* \left(\ch \left({\bf R}(p'' \times p')_*
(G'' \otimes {\cal R} \otimes {G'}^{\vee})\right) \right)\cdot\\
& \quad \quad \pi_{X' \times X}^* \left( \ch \left({\bf R}(p' \times p)_*
(G' \otimes {\cal Q} \otimes {G}^{\vee})\right) \right)
\pi_{X'}^*\left(\frac{\td_{X'}}
{\ch({\bf R}p'_*({G'}^{\vee} \otimes G'))}\right)
\\
=& \ch \left({\bf R}(p'' \times p' \times p)_*
(G'' \otimes {\cal R} \otimes {\cal Q} \otimes G^{\vee})\right)
\pi_{X'}^*(\td_{X'}).
\end{split}
\end{equation}
Since
\begin{equation}
\begin{split}
&\pi_{X'' \times X*}\left(\ch \left({\bf R}(p'' \times p' \times p)_*
(G'' \otimes {\cal R} \otimes {\cal Q} \otimes G^{\vee})\right)
\pi_{X'}^*(\td_{X'})\right)\\
=&\ch\left({\bf R}\pi_{X'' \times X*}\left({\bf R}(p'' \times p' \times p)_*
(G'' \otimes {\cal R} \otimes {\cal Q} \otimes G^{\vee})\right)\right)\\
=& \ch({\bf R}(p'' \times p)_* \circ {\bf R}\tilde{\pi}_{Y'' \times Y*}
(G'' \otimes {\cal R} \otimes {\cal Q} \otimes G^{\vee}))\\
=& \ch({\bf R}(p'' \times p)_*(G'' \otimes {\cal S} \otimes G^{\vee})),
\end{split}
\end{equation}
we get

\begin{equation}
\begin{split}
\Psi_{X'' \to X}^{\cal S}(z)=&
\pi_{X*}\left( \pi_{X''}^*(z)  
\ch({\bf R}(p'' \times p)_*({G''}  \otimes {\cal S} \otimes  G^{\vee}))
\frac{\sqrt{\td_{X''}}}{\sqrt{\ch({\bf R}p''_*({G''}^{\vee} \otimes G''))}}
\frac{\sqrt{\td_X}}{\sqrt{\ch({\bf R}p_*(G^{\vee} \otimes G))}}\right)\\
=&
\Psi_{X' \to X}^{\cal Q} \circ \Psi_{X'' \to X'}^{\cal R}(z).
\end{split}
\end{equation}
\end{proof}

\begin{lem}
Assume that the canonical bundles $K_X, K_{X'}$ are trivial.
Then
\begin{equation}
\langle x,\Psi_{X' \to X}^{\cal Q}(y)\rangle 
=\langle \Psi_{X \to X'}^{{\cal Q}^{\vee}}(x), y \rangle,\;\;
x \in H^*(X,{\Bbb Q}),\; y \in H^*(X',{\Bbb Q}),
\end{equation}
where $\langle\;\;,\;\; \rangle$ is the Mukai pairing.
\end{lem}

\begin{proof}
\begin{equation}
\begin{split}
&\langle x,\Psi_{X' \to X}^{\cal Q}(y) \rangle \\
=&-\int_{X} x \Psi_{X' \to X}^{\cal Q}(y)^{\vee}\\
=&-\int_{X' \times X} \pi_X^*(x) 
\left(\pi_{X'}^*(y)
\frac{\sqrt{\td_{X'}}}{\sqrt{\ch({\bf R}p'_*({G'}^{\vee} \otimes G'))}} 
\frac{\sqrt{\td_X}}{\sqrt{\ch({\bf R}p_*(G^{\vee} \otimes G))}}
\ch({\bf R}(p' \times p)_*({G'} \otimes {\cal Q} \otimes  G^{\vee}))
\right)^{\vee}
\\
=&-\int_{X' \times X} 
\left(\frac{\sqrt{\td_{X'}}}{\sqrt{\ch({\bf R}p'_*({G'}^{\vee} \otimes G'))}} 
\frac{\sqrt{\td_X}}{\sqrt{\ch({\bf R}p_*(G^{\vee} \otimes G))}}
\ch({\bf R}(p' \times p)_*({G'}^{\vee}\otimes {\cal Q}^{\vee} \otimes  G ))
\pi_X^*(x)\right) \pi_{X'}^*(y^{\vee}) \\
=&-\int_{X'} \Psi_{X \to X'}^{{\cal Q}^{\vee}}(x) y^{\vee}\\
=&\langle \Psi_{X \to X'}^{{\cal Q}^{\vee}}(x),y \rangle.
\end{split}
\end{equation}
\end{proof}

\subsection{Fourier-Mukai transform induced by stable twisted sheaves}

Let $p:Y \to X$ be a projective bundle over 
an abelian surface or a $K3$ surface.
Let $G$ be a locally free $Y$-sheaf.
%Let $\alpha=\{\alpha_{ijk}\}$ be a 2-cocycle which is torsion. 
Assume that $X':=\overline{M}^{Y,G}_H(v)$ is a surface and
consists of stable sheaves.
We set $Y':=\widetilde{\overline{M}^{Y,G}_H(v)}$.
Let ${\cal E}$ be the family on $Y' \times Y$.
\begin{NB}
Hence there is an open covering $M_H(v)=\cup_j V_j$
and local universal families ${\cal E}_j$ on $V_j \times Y$.
By using an open covering of $X$,
${\cal E}_j$ is a collection of
${\cal E}_{(i,j)} \in \Coh(V_j \times U_i)$
and homomorphisms $\phi_{i_2,i_1}(j):
{\cal E}_{(i_2,j)|V_j \times (U_{i_2} \cap U_{i_1})} \to
 {\cal E}_{(i_1,j)|V_j \times (U_{i_1} \cap U_{i_2})}$
such that $\phi_{i_2,i_1}(j) \circ \phi_{i_3,i_2}(j) 
\circ \phi_{i_1,i_3}(j)=\alpha_{i_1 i_2 i_3}$.
Since $\Hom_{\pi_{V_i}}({\cal E}_i,{\cal E}_i) \cong {\cal O}_{V_i}$,
there is a 2-cocycle $\beta:=\{\beta_{ijk}|
\beta_{ijk} \in H^0(V_i \cap V_j \cap V_k,{\cal O}_{X'}^{\times}) \}$
such that ${\cal E}_i$ is a $\pi_{X'}^*(\beta)$-twisted sheaf on 
$X' \times Y$. That is
${\cal E}:=\{(V_i \times X,{\cal E}_i)\}$ 
is a $\pi_{X'}^*(\beta)\pi_X^*(\alpha)$-twisted sheaf on 
$X' \times X$.
A homomorphism $\varphi_{j_2,j_1}:{\cal E}_{j_2} \to {\cal E}_{j_1}$
is a collection of $\varphi_{j_2,j_1}(i):{\cal E}_{(i,j_2)} \to
{\cal E}_{(i,j_1)}$ such that the following diagrams are commutative:
\begin{equation}
\begin{CD}
{\cal E}_{(i_2,j_2)|(V_{j_1} \cap V_{j_2}) \times (U_{i_2} \cap U_{i_1})}
@>{\phi_{i_2,i_1}(j_2)}>>
{\cal E}_{(i_1,j_2)|(V_{j_1} \cap V_{j_2}) \times (U_{i_2} \cap U_{i_1})}\\
@V{\varphi_{j_2,j_1}(i_2)}VV @VV{\varphi_{j_2,j_1}(i_1)}V\\
{\cal E}_{(i_2,j_1)|(V_{j_1} \cap V_{j_2}) \times (U_{i_2} \cap U_{i_1})}
@>>{\phi_{i_2,i_1}(j_1)}>
{\cal E}_{(i_1,j_1)|(V_{j_1} \cap V_{j_2}) \times (U_{i_2} \cap U_{i_1})}.
\end{CD}
\end{equation}
We set $\psi_{(i_2,j_2),(i_1,j_1)}:=
\varphi_{j_2,j_1}(i_1) \circ \phi_{i_2,i_1}(j_2)$.
Then $\psi_{(i_2,j_2),(i_1,j_1)} \circ 
\psi_{(i_3,j_3),(i_2,j_2)}\circ \psi_{(i_1,j_1),(i_3,j_3)}=
\alpha_{i_1 i_2 i_3}\beta_{j_1 j_2 j_3}$.
 \end{NB}

%We set $w(Y):=[\xi \mod r]$, $\xi \in H^2(X,{\Bbb Z})$ and
%$w(Y'):=[\xi' \mod r]$, $\xi' \in H^2(X',{\Bbb Z})$.

We consider integral functors
\begin{equation}
\begin{matrix}
\Phi_{X' \to X}^{\cal E}: &{\bf D}(X',Y')& \to &{\bf D}(X,Y)\\
& x & \mapsto & 
{\bf R}\widetilde{\pi}_{Y*}
({\cal E} \otimes \widetilde{\pi}_{Y'}^*(x)),
\end{matrix}
\end{equation}

\begin{equation}
\begin{matrix}
\Phi_{X \to X'}^{{\cal E}^{\vee}}[2]:
&{\bf D}(X,X)& \to &{\bf D}(X',Y')\\
& y & \mapsto & 
{\bf R}\widetilde{\pi}_{Y'*}
({\cal E}^{\vee} \otimes \widetilde{\pi}_{Y}^*(y)[2]).
\end{matrix}
\end{equation}

\begin{rem}
Let ${\cal L}({p'}^*({\alpha}^{-1}))$ and 
${\cal L}(p^*(\alpha^{-1}))$ be twisted line bundles
on $Y'$ and $Y$ respectively in \eqref{eq:Y=X}. 
%which give equivalences
%$\Lambda^{{\cal L}({p'}^*({\alpha'}^{-1}))}:
%{\bf D}(X',Y') \cong {\bf D}(X',\alpha')$ and
%$\Lambda^{{\cal L}(p^*({\alpha}^{-1}))}:
%{\bf D}(X,Y) \cong {\bf D}(X,\alpha)$ in \eqref{eq:Y=X}.
Then $\Lambda^{{\cal L}(p^*({\alpha}^{-1}))} \circ
\Phi_{X' \to X}^{\cal E} \circ
(\Lambda^{{\cal L}({p'}^*({\alpha'}^{-1}))})^{-1}:
{\bf D}(X',\alpha') \to {\bf D}(X,\alpha)$ is an integral functor
with the kernel
${\bf R}(p' \times p)_*({\cal L}({p'}^*({\alpha'}^{-1})) \otimes 
{\cal E} \otimes
{\cal L}(p^*({\alpha}^{-1}))^{\vee}) \in
{\bf D}(X' \times X,{\alpha'}^{-1} \times \alpha)$.
\end{rem}

C\u ald\u araru \cite{C:2} developed a theory of derived category
of twisted sheaves. 
In particular, Grothendieck-Serre duality holds.
Then we see that
$\Phi_{X \to X'}^{{\cal E}^{\vee}}[2]$ is the adjoint of 
$\Phi_{X' \to X}^{{\cal E}}$.
As in the usual Fourier-Mukai functor, we see that
the following theorem holds
(see \cite{Br:1}, \cite{C:1}). 
\begin{thm}\label{thm:Fourier-Mukai}
$\Phi_{X \to X'}^{{\cal E}^{\vee}}[2] \circ 
\Phi_{X' \to X}^{\cal E} \cong 1$ and 
$\Phi_{X' \to X}^{\cal E} \circ \Phi_{X \to X'}^{{\cal E}^{\vee}}[2] \cong 1$.
Thus
$\Phi_{X' \to X}^{\cal E}$ is an equivalence.
\end{thm}

Then we have the following which also follows from a more
general statement \cite[Thm. 0.4]{H-S:2}.

\begin{cor}
$\Psi_{X' \to X}^{\cal E}$ induces an isometry of the 
Hodge structures:
\begin{equation}
(H^*(X',{\Bbb Z}),\langle\;\;,\;\; \rangle,-\frac{\xi'}{r})
\cong 
(H^*(X,{\Bbb Z}),\langle\;\;,\;\; \rangle,-\frac{\xi}{r}).
\end{equation}
\end{cor}

\begin{proof}
Obviously
$\Psi_{X' \to X}^{\cal E}$ induces an isometry of
the Hodge structures over ${\Bbb Q}$.
If $X$ is a $K3$ surface such that $w(Y) \in \NS(X) \otimes \mu_r$
and $X'$ is a fine moduli space, then
$\Psi_{X' \to X}^{\cal E}$ is defined over ${\Bbb Z}$.
For a general case, we use the deformation arguments.
\end{proof} 
We also have the following which is used in \cite{Y:12}.

\begin{cor}
Assume that $X'$ consists of locally free $Y$-sheaves.
Then ${\cal E}^{\vee}_{|Y' \times \{y \}}$, $y \in Y$
is a simple $Y'$-sheaf.
If $\NS(X) \cong {\Bbb Z}H$, then 
${\cal E}^{\vee}_{|Y' \times \{y \}}$, $y \in Y$
is a stable $Y'$-sheaf.
\end{cor} 

\begin{proof}
Since $\Phi_{X \to X'}^{{\cal E}^{\vee}}[2]$ is an equivalence,
$\Phi_{X \to X'}^{{\cal E}^{\vee}}({\cal O}_{p^{-1}(p(y))}(1))
={\cal E}^{\vee}_{|Y' \times \{y \}}$ 
is a simple $Y'$-sheaf.
If $\NS(X) \cong {\Bbb Z}$, then Proposition \ref{prop:simple}
implies the stability of ${\cal E}^{\vee}_{|Y' \times \{y \}}$.
\end{proof}

\begin{NB}
Let $p':Y' \to X$ be a projective bundle.
Let ${\cal E}^1,{\cal E}^2$ be pull-backs of
${\cal E}$ by $Y' \times Y \times Y' \to  Y' \times Y$
and $Y' \times Y \times Y' \to  Y \times Y'$
respectively.
 
$(p' \times p')_*
({\bf R}\Hom_{p_{Y' \times Y'}}({\cal E}^2 \otimes {G'}^{\vee},
{\cal E}^1\otimes {G'}^{\vee}))[2]
\cong {\cal O}_{\Delta_X} \otimes p'_*({G'}^{\vee} \otimes G')$.
${\cal O}_{Y_i'}(\lambda'_i) \boxtimes_{{\cal O}_{V_i}} 
{\cal O}_{Y_i'}(-\lambda'_i)$ with a transition function
$\phi_{ji} \boxtimes_{{\cal O}_{V_i \cap V_j}} {}^t\phi_{ji}^{-1}$
defines a line bundle ${\cal L}$ on $Y' \times_{X'} Y'$. 

${\cal E}_{|Y'_i \times Y} \cong (p' \times 1)^*({\cal E}_i)
\otimes {\cal O}_{Y'_i}(\lambda'_i)$.
${\bf R}\Hom_{p_{V_i' \times V_i'}}({\cal E}_i^2,
{\cal E}_i^1)[2]
\cong {\cal O}_{\Delta_{V_i}}$ and it defines a structure sheaf of the
diagonal. 
$G'_{|Y_i'}={p'}^*(G'_i)(\lambda'_i)$.
Hence ${\cal L} \otimes {G'}^{\vee} \boxtimes_{{\cal O}_X} {G'}_{|Y_i}
= {G'_i}^{\vee} \boxtimes_{{\cal O}_{V_i}} {G'_i}$.
Hence $(p' \times p')_*
({\cal L} \otimes {G'}^{\vee} \boxtimes_{{\cal O}_X} {G'})={\cal O}_{\Delta_X}
\otimes p'_*({G'}^{\vee} \otimes G')$.
\end{NB}

\vspace{1pc}

{\it Acknowledgement.}
First of all, I would like to thank
Daniel Huybrechts and Paolo Stellari.
They proved C\u ald\u araru's conjecture.
Moreover Huybrechts gave me many
valuable suggestions on this paper.
I would also like to thank Eyal Markman and Shigeru Mukai
for valuable discussions on the twisted sheaves and their moduli spaces.
Thanks also to
Max Lieblich for explaining the relation of our moduli spaces
with Simpson's moduli spaces of modules over the 
Azumaya algebra.


\begin{thebibliography}{[G-H]}
\bibitem[Br]{Br:1}
Bridgeland, T.,
{\it Equivalences of triangulated categories and Fourier-Mukai
transforms,} 
Bull. London Math. Soc. {\bf 31} (1999), 25--34,
math.AG/9809114

\bibitem[C1]{C:1}
C\u ald\u araru, A.,
{\it Nonfine moduli spaces of sheaves on $K3$ surfaces,}
Int. Math. Res. Not. 2002, no. 20, 1027--1056
\bibitem[C2]{C:2}
C\u ald\u araru, A.,
{\it Derived categories of twisted sheaves on Calabi-Yau manifolds,}
Ph.D. thesis, Cornell University (2000)

\bibitem[D]{D:1}
De Jong, A. J.,
{\it The period-index problem for the Brauer group of an algebraic surface,}
Duke Math. J. {\bf 123} (2004), 71--94
\bibitem[Ho-St]{Ho-St:1}
Hoffmann, N., Stuhler, U.,
{\it Moduli schemes of rank one Azumaya modules,}
math.AG/0411094
\bibitem[H-L]{H-L:1}
Huybrechts, D., Lehn, M., 
{\it The geometry of moduli spaces of sheaves,}
Aspects of Mathematics, E31. Friedr. Vieweg \& Sohn, Braunschweig, 1997
\bibitem[H-Sc]{H-S:1}
Huybrechts, D., Schr\"{o}er, S.,
{\it The Brauer group of analytic $K3$ surfaces,}
Int. Math. Res. Not. (2003), no. 50, 2687--2698
\bibitem[H-St]{H-S:2}
Huybrechts, D., Stellari, P.,
{\it Equivalences of twisted $K3$ surfaces,}
math.AG/0409030
\bibitem[H-St2]{H-St:2}
Huybrechts, D., Stellari, P.,
{\it Proof of C\u ald\u araru's conjecture.
An appendix to a paper by Yoshioka,}
math.AG/0411541
\bibitem[L]{Langer:1}
Langer A.,
{\it Moduli spaces of sheaves in mixed characteristic,}
Duke Math. J., {\bf 124} (2004), 571--586
\bibitem[Li]{Li:1}
Lieblich, M.,
{\it Moduli of twisted sheaves,}
math.AG/0411337
 \bibitem[Mu1]{Mu:2}
Mukai, S.,
{\it Duality between $D(X)$ and $D(\hat{X})$ with its application
to Picard sheaves,}
Nagoya Math. J., {\bf 81} (1981), 153--175
\bibitem[Mu2]{Mu:3}
Mukai, S.,
{\it Symplectic structure of the moduli space of sheaves on an 
abelian or K3 surface,}
Invent. math. {\bf 77}
(1984), 101--116
\bibitem[Mu3]{Mu:4}
Mukai, S.,
{\it On the moduli space of bundles on K3 surfaces I,}
Vector bundles on Algebraic Varieties, Oxford, 1987, 341--413
\bibitem[Mu4]{Mu:5}
Mukai, S.,
{\it Vector bundles on a $K3$ surface,}
Proceedings of the International Congress of Mathematicians, 
Vol. II (Beijing, 2002), 495--502, Higher Ed. Press, Beijing, 2002. 
\bibitem[Or]{Or:1}
Orlov, D.,
{\it Equivalences of derived categories and K3 surfaces,}
J. Math. Sci. (NY), {\bf 84} (1997) 1361--1381 
\bibitem[S]{S:1}
Simpson, C.,
{\it Moduli of representations of the fundamental group
of a smooth projective variety I,}
Publ. Math. I.H.E.S. {\bf 79} (1994), 47--129                                                
\bibitem[Y1]{Y:7}
Yoshioka, K.,
{\it Moduli spaces of stable sheaves on abelian surfaces,}
math.AG/0009001, Math. Ann. {\bf 321} (2001), 817--884
\bibitem[Y2]{Y:9}
Yoshioka, K.,
{\it Twisted stability and Fourier-Mukai transform I,}
Compositio Math. {\bf 138} (2003), 261--288
\bibitem[Y3]{Y:11}
Yoshioka, K.,
{\it Twisted stability and Fourier-Mukai transform II,}
Manuscripta Math. {\bf 110} (2003), 433--465
\bibitem[Y4]{Y:12}
Yoshioka, K.,
{\it Stability and the Fourier-Mukai transform II,}
preprint
(sections 3, 4 of math.AG/0112267)


\begin{NB}
\bibitem[Ma]{Mar:4}
Maruyama, M.,
{\it Construction of moduli spaces of stable sheaves 
via Simpson's idea,}
 Moduli of vector bundles (Sanda, 1994; Kyoto, 1994), 147--187, 
Lecture Notes in Pure and Appl. Math., 179, Dekker, New York, 1996 
\bibitem[Mum]{Mum:1}
Mumford, D.,
{\it Lectures on curves on an algebraic surface,}
Annals of Mathematics Studies, No. 59 Princeton University Press, 
Princeton, N.J. 1966 
\end{NB}

\end{thebibliography}
\end{document}